\documentclass[11pt]{article}
\usepackage{amssymb}
\usepackage{amsmath}

\newtheorem{theo}{Theorem}[section]

\newtheorem{prop}[theo]{Proposition}
\newtheorem{lemm}[theo]{Lemma}
\newtheorem{coro}[theo]{Corollary}
\newtheorem{rema}[theo]{Remark}

\newtheorem{conj}[theo]{Conjecture}

\newcommand{\cqfd}
{%
\mbox{}%
\nolinebreak%
\hfill%
\rule{2mm}{2mm}%
\medbreak%
\par%
}
\newfont{\gothic}{eufb10}
\date{\empty}
\begin{document}
\title{Infinitesimal invariants for cycles modulo  algebraic equivalence and $1$-cycles on Jacobians}
\author{Claire Voisin\\CNRS and \'{E}cole Polytechnique} \maketitle \setcounter{section}{-1}
\begin{abstract} We construct an infinitesimal invariant for
cycles in a family with cohomology class in the total space lying in
a given level of the Leray filtration. This infinitesimal invariant
detects cycles modulo algebraic equivalence in the fibers. We apply
this construction to
 the Ikeda family, which  gives  optimal results for the Beauville decomposition of the
$1$-cycle of a very general plane curve in its Jacobian.
\end{abstract}
\section{Introduction}
\setcounter{equation}{0} This paper is devoted  first of all to the
construction of  infinitesimal invariants for (families of) cycles
modulo algebraic equivalence, similar to the infinitesimal
invariants used in \cite{voisinvar}  for the study of families of
cycles modulo rational equivalence, and secondly to a geometric
application of these infinitesimal invariants to the study of the
length of the Beauville decomposition of the $1$-cycle of a curve in
its Jacobian.

Infinitesimal invariants appear in the following situation: Let
$\pi:\mathcal{Y} \rightarrow B$ be a smooth projective morphism, and
let $\mathcal{Z}$ be a codimension $n$ cycle in $\mathcal{Y}$. We
assume that the cohomology class $[\mathcal{Z}]\in
H^{2n}(\mathcal{Y},\mathbb{Q})$ belongs to the $s$-th level of the
Leray filtration on  $H^{2n}(\mathcal{Y},\mathbb{Q})$. Then we
produce for any point $b\in B$ an infinitesimal invariant
$\delta[\mathcal{Z}]_{b,alg}$ which depends only on  the first order
neighborhood $\mathcal{Z}_{b,2}=\mathcal{Z}_{\mid
\mathcal{Y}_{b,2}}$ of the restricted cycle $\mathcal{Z}_{b}$ of the
fibre $\mathcal{X}_b$ (hence the name ``infinitesimal invariant'').
More precisely, $\mathcal{Z}$ has a Dolbeault cohomology class
$[\mathcal{Z}]^{n,n}\in H^n(\mathcal{Y},\Omega_{\mathcal{X}}^n)$ and
this infinitesimal invariant will depend only on the restricted
class
$$[\mathcal{Z}]^{n,n}_{\mid \mathcal{Y}_b}\in
H^n(\mathcal{Y}_b,{\Omega_{\mathcal{Y}}^n}_{\mid \mathcal{Y}_b}).$$
This infinitesimal invariant  has the property that it allows to
decide whether the restricted cycle $\mathcal{Z}_{b}$ is
algebraically equivalent to $0$ at the very general point $b\in B$.
More precisely, we will prove the following result:
\begin{theo} \label{theoinfiintro} In the situation above, assuming $s\geq2$ and $[\mathcal{Z}]\in L^sH^{2n}(\mathcal{Y},\mathbb{Q})$, there is
an infinitesimal invariant
$$\delta[\mathcal{Z}]_{alg}\in \frac{\text{Ker}\,\overline{\nabla}_{s,n-s,n}}{\Omega_{B}\wedge
\text{Ker}\,\overline{\nabla}_{s-1,n-s,s}},$$ with restriction
$$\delta[\mathcal{Z}]_{b,alg}\in \frac{\text{Ker}\,\overline{\nabla}_{s,n-s,n,b}}{\Omega_{B,b}\wedge \text{Ker}\,\overline{\nabla}_{s-1,n-s,s,b}} ,$$
which is $0$ at a general point $b\in B$ if for general $b\in B$,
  $\mathcal{Z}_{b}$ is algebraically equivalent to $0$  up to torsion in the fiber
$\mathcal{Y}_b$.
\end{theo}
The maps $$\overline{\nabla}_{s,n-s,n}:\Omega_{B}^{s}\otimes
\mathcal{H}^{n-s,n}\rightarrow \Omega_{B}^{s+1}\otimes
\mathcal{H}^{n-s-1,n+1}$$
 are induced by the infinitesimal
variation  of Hodge structure  on the cohomology of the fibers of
$\pi$; they are defined using the Gauss-Manin connection and the
Hodge filtration on Hodge bundles associated to the family
$\mathcal{Y}\rightarrow B$ (cf. \cite{griffithsIVHS},
\cite[5.1.2]{voisinbook}).

Since the set of points $b\in B$ such that a multiple of the
restricted cycle $\mathcal{Z}_b$ is algebraically equivalent to $0$
is a countable union of closed algebraic subsets of $B$, we can
rephrase Theorem \ref{theoinfiintro} as follows: \begin{coro}
\label{coroinfiintro} There is a Zariski open set $U\subset B$ such
that, if $\delta[\mathcal{Z}]_{b,alg}\not=0$ for some $b\in U$, then
$[\mathcal{Z}]_{b'}$ is not algebraically equivalent to $0$ up to
torsion for very general $b'\in B$. \end{coro}

In fact we can take for $U$ the Zariski open set of $B$ where the
sheaf
$\frac{\text{Ker}\,\overline{\nabla}_{s,n-s,n}}{\Omega_{B,b}\wedge
\text{Ker}\,\overline{\nabla}_{s-1,n-s,s}}$ is locally free; the non
vanishing of $\delta[\mathcal{Z}]_{b,alg}$ for some $b\in U$ then
guarantees the non vanishing of $\delta[\mathcal{Z}]_{b,alg}$ for
general $b\in B$.

 The main difference
between the present construction and the one in \cite{voisinvar}
lies in the fact that the infinitesimal invariant introduced here
detects cycles modulo algebraic equivalence, while in our previous
work only cycles modulo rational equivalence (a much bigger group)
could be detected. There is a serious difficulty here, since Lemma
\ref{lemmarat}, which says that up to shrinking the base $B$, the
cycle $\mathcal{Z}$ as above is cohomologous to $0$ modulo torsion
on the total space $\mathcal{Y}$ if the restricted cycles
$\mathcal{Z}_t$ are rationally equivalent to $0$, obviously does not
hold when replacing rational equivalence by algebraic equivalence.

Our infinitesimal invariants are not present, at least explicitly,
in the work of Nori \cite{nori} and Fakhruddin \cite{fakhruddin},
but they are in fact  hidden in the infinitesimal variations of
Hodge structures arguments used in their computations of the Leray
filtration (we refer to \cite{voisinicm}, \cite{voisinbook} for the
relation between infinitesimal invariants and Leray filtration).
 The work of Ikeda \cite{ikeda} involves similar infinitesimal invariants, which lie however in
  a bigger cohomology group, and allow to detect the nontriviality of cycles
 modulo a relation which is much finer than algebraic equivalence.
  His computations have been extended by
 Pirola and Rizzi in \cite{pirola}.

Our main application concerns the lenght of the Beauville
decomposition of the canonical $1$-cycle of a generic Jacobian. To
avoid heavy notation, we denote by $\text{CH}^i(Y)$ the Chow groups
of a variety $Y$  with rational coefficients. Recall that Beauville
defined in \cite{beauvilledecom} a canonical decomposition on the
Chow groups $\text{CH}^i(A)$ of an abelian variety, as a direct sum
of eigenspaces for the actions of the homotheties:
$$\text{CH}^i(A)=\oplus_s\text{CH}^i(A)_s,$$ where
$$z\in \text{CH}^i(A)_s\Leftrightarrow \mu_m^*z=m^{2i-s}z,\,\forall m\in \mathbb{Z}.$$
Here for $m\in \mathbb{Z}$ we denote by $\mu_m$ the multiplication
by $m$ on $A$. This decomposition works as well for cycles modulo
algebraic equivalence. It should be understood as  a canonical
splitting of the conjectural Bloch-Beilinson filtration, which is
supposed to be  a decreasing filtration $F^l\text{CH}^i(A)$ on
$\text{CH}^i(A)$, whose graded pieces $Gr_F^l\text{CH}^i(A)$ are
governed by the Hodge structures on the cohomology groups
$H^{2i-l}(A,\mathbb{Q})$, and which finishes, more precisely
$F^{i+1}\text{CH}^i(A)=0$.

  It is not known however if the
filtration defined by
\begin{eqnarray}
\label{eqbeaudec}F^l\text{CH}^i(A)=\oplus_{s\geq l}\text{CH}^i(A)_s
\end{eqnarray}
satisfies the axioms of a Bloch-Beilinson filtration, and in
particular if
$$F^1\text{CH}^i(A)=\text{CH}^i(A)_{hom}.$$
Here the inclusion $F^1\text{CH}^i(A)\subset \text{CH}^i(A)_{hom}$
is easy to see : for $s>0$, $\text{CH}^i(A)_s$ is made of cycles
homologous to $0$, because $\mu_n^*$ acts as multiplication by
$n^{2i}$ on $H^{2i}(A,\mathbb{Q})$, but the reverse inclusion, which
is equivalent to saying that the cycle class map is injective on
$\text{CH}^i(A)_0$, is unknown.

In this paper, we consider the  induced decomposition on the  groups
$\text{CH}^i(A)/ alg$ of cycles with rational coefficients modulo
algebraic equivalence. As mentioned above, for $s>0$,
$\text{CH}^i(A)_s\subset \text{CH}^i(A)_{hom}$, so that the $s$-th
piece $\text{CH}^i(A)_s/ alg$ is in fact contained in the Griffiths
group
$$ \text{CH}^i(A)_{hom}/ alg=\text{Griff}^i(A)\otimes\mathbb{Q}.$$
We will  avoid using the notation $\text{Griff}^i(A)_s$ since it is
used by Ikeda with a different meaning (see below). In the case
where $A=JC$ is the Jacobian of a curve of genus $g$, there is a
natural $1$-cycle (codimension $g-1$ cycle) $Z$ on $JC$ defined as
the image of $C$ under the Abel  map, with respect to any chosen
point or $0$-cycle of degree $1$ on $C$. Being defined up to
translation, it is well-defined modulo algebraic equivalence. We
consider the Beauville decomposition
$$Z=\sum_s Z_s$$
of this cycle, with $Z_s\in \text{CH}^{g-1}(JC)_s/ alg$. The study
of the components $Z_s$ has a long history. The first result is the
one by Ceresa \cite{ceresa}, saying that $Z_1$ is non zero in
$\text{CH}^{g-1}(JC)_1/ alg$ for $g\geq 3$ and $C$ very general in
modulus. Ceresa uses
  the Griffiths criterion
\cite{griffiths} saying that if a cycle  is algebraically equivalent
to $0$, its Abel-Jacobi invariant belongs to the maximal abelian
subvariety of the intermediate Jacobian, together with a
degeneration argument for the study of this Abel-Jacobi invariant.
 Using Griffiths' infinitesimal invariant \cite{griffithsIVHS}, an infinitesimal argument    can be applied  as well (cf. \cite{collinopirola})
 to prove the nontriviality of the primitive (or transcendental) part of the normal function given by the Abel-Jacobi
 image of the cycle $Z_1$.

In the opposite direction, there is  the following vanishing result
due to  Colombo and Van Geemen \cite{geemen}:
\begin{theo} \label{theocvg} Assume a smooth curve $C$ has a $g_{d}^1$, that is a morphism of degree $d$ to $\mathbb{P}^1$. Then
$$Z_s=0\,\,{\rm in}\,\,\text{CH}^{g-1}(JC)_s/ alg,\,\forall s\geq d-1.$$
\end{theo}
Further results concerning the so-called ``tautological ring''
introduced by Beauville in \cite{beauvillejac} and generated under
Pontryagin product by the $Z_s$ have been obtained by Herbaut
\cite{herbaut}, Van der Geer-Kouvidakis \cite{vandergeer},
Polishchuk and Moonen \cite{moonenpoli}, \cite{moonen}. The results
of Herbaut (reproved by Van der Geer-Kouvidakis) generalize the
Colombo-Van Geemen result to the case where the curve admits a
$g_d^r$, $r\geq 2$. This however does not give more vanishing for
the individual components $Z_s$, but rather polynomial relations
 between them in the tautological ring (here the ring structure considered is given by Pontryagin product of cycles).
The results by Moonen and Polishchuk are devoted to relations in the
tautological ring of a general curve, and the paper
\cite{beauvillejac} by Beauville describes a set of generators for
this ring.

In the paper \cite{geemen}, Colombo and Van Geemen asked for the
optimality of their result, for a very general curve. (As noticed by
Ben Moonen, it is not expected that their result is optimal for any
curve, because for curves defined over a number field, the Beilinson
conjecture \cite{beilinson} predicts
$$Z_s=0,\,\forall s\geq 2$$
while the generic curve defined over a number field and of genus $\geq 5$ has gonality $\geq 4$.)

Concretely, the gonality $\text{gon}(C)$ of a generic curve $C$ of
genus $g$, that is the minimal degree of a $g_d^1$ on $C$ is
computed by Brill-Noether theory (cf. \cite{acgh}) which gives:
$$ g=2k+1 \,\,{\rm or}\,\,g=2k+2\Rightarrow \text{gon}(C)=k+2\,\,{\rm for\,\, generic}\,\,C.$$
Thus the question asked by Colombo and Van Geemen ican be stated in
the following form :
\begin{conj}\label{conj} Let $C$ be a very general curve of genus $g\geq 2k+1$. Then
$$Z_k\not=0\,\,{\rm in}\,\, \text{CH}^{g-1}(JC)_{k}/ alg.$$

\end{conj}
Note that by a degeneration argument as in Ceresa's paper
\cite{ceresa}, we can reduce to the case where $g=2k+1$. Some
results have been obtained in this direction. First of all
Fakhruddin \cite{fakhruddin} shows the non vanishing of $Z_2$ for
$C$ a general curve of genus $\geq 11$ (the expected bound being
$g\geq 5$). In a different direction,  Ikeda obtains a non vanishing
result for an analogue of the cycles $Z_s,\,s\leq d-2$ when the
curve is a very general plane curve of degree $d$ (thus of gonality
$d-1$). Colombo-van Geemen's Theorem \ref{theocvg} says that for
such a curve, the cycle $Z_{d-2}$ is zero in
$\text{CH}^{g-1}(JC)_{d-2}/ alg$. This last result is also possibly
optimal since by  Herbaut's above mentioned results in
\cite{herbaut}, we know that for plane curves, there are extra
polynomial relations between the various $Z_k$, but a priori no
further vanishing of a given component $Z_k$ is expected except for
the vanishing of $Z_{d-2}$ given by Theorem \ref{theocvg}. This
leads to the following conjecture:
\begin{conj}\label{conj2} Let $C$  be a very general plane curve of
degree $d$ (so $g=g(C) $ is given by $2g-2=d(d-3)$). Then the cycle
$Z_{d-3}$ is $\not=0\,\,{\rm in}\,\, \text{CH}^{g-1}(JC)_{d-3}/
alg$.
\end{conj}
Ikeda addresses this conjecture but not for cycles modulo algebraic
equivalence.
 The quotient of the group of
cycles considered by Ikeda is  an a priori  much bigger group than
the group of cycles modulo algebraic equivalence. These groups
called ``higher Griffiths groups'' were defined by Shuji Saito
\cite{saito} and involve his construction of a Bloch-Beilinson
filtration $F^i_S\text{CH}^k(Y)$ \cite{ssaito}. The ``higher
Griffiths groups'' of Saito are then defined as the quotient
$$F^i_S\text{CH}^k(Y)/F^{i+1}_S\text{CH}^k(Y)+<\Gamma_*F^i_SCH_0(Z)>,$$ where $\Gamma$ runs over the
set of correspondences of codimension $k$ between any smooth
projective variety $Z$ and $Y$. Of course when $i =1$, the subgroup
$<\Gamma_*F^1_S\text{CH}_0(Z)>$ consists exactly of cycles
algebraically equivalent to $0$, but for   $i\geq 1$, the group
$<\Gamma_*F^i_SCH_0(Z)>\subset F^i_S\text{CH}^k(Y)$ may be very
small, and different from $F^i_S\text{CH}^k(Y)\cap
\text{CH}^k(Y)_{alg}$. We will make a comparison between the Ikeda
infinitesimal invariants and ours in section \ref{secdiscussion},
and we will see that, due to the vanishing of canonical syzygies for
generic curves \cite{voisinsyz},  for the case of a general curve of
genus $2k+1$, his infinitesimal invariants for the cycle $Z_s$,
$s\leq k$, highly depend on the choice of embedding of $C$ into its
Jacobian (this embedding, hence also the cycles $Z_s$, is defined up
to translation).

 Our main application of the infinitesimal invariant constructed in
Section
  \ref{sec1} is the proof of conjecture \ref{conj2}.
\begin{theo}\label{main} Let $C$ be a very general plane curve of degree $d$. Then
$$Z_{d-3}\not=0\,\,{\rm in}\,\, \text{CH}^{g-1}(JC)_{d-3}/ alg.$$

\end{theo}
Theorem \ref{main} implies in particular the following result, which
is certainly not optimal in view of Conjecture \ref{conj}, but
answers partially the question asked in \cite{beauvillejac}:
\begin{coro} Let $C$ be a very general curve of genus $g> \frac{k(k+3 )}{2}$. Then the
cycle $Z$ satisfies $Z_k\not=0$ in $\text{CH}^{g-1}(JC)_k$.
\end{coro}
{\bf Proof.} Indeed, we first reduce  by degeneration  the problem
to the case where $g= \frac{k(k+3 )}{2}+1$; then we specialize $C$
to a very general plane curve of degree $k+3$ and apply Theorem
\ref{main}. Since the specialized cycle $Z_k$ is not algebraically
equivalent to $0$, the same is true at the very general point.
\cqfd
 In order to prove Theorem \ref{main}, we compute the
infinitesimal invariant of Theorem \ref{theoinfiintro} on the Ikeda
family of plane curves of degree $d$ (cf. \cite{ikeda}). This family
has a very special variation of Hodge structure, which makes the
explicit computation quite easy.

 \section{Infinitesimal invariants for cycles modulo algebraic equivalence \label{sec1}}
 \subsection{Infinitesimal invariants for cycles modulo rational rational equivalence}
 \label{subsecrat}
 In this subsection,  we recall for the convenience of the reader the construction of  infinitesimal invariants for families of cycles modulo
 rational equivalence.
 We will explain in next subsection \ref{secalgeqinf} the modifications
 needed to get infinitesimal invariants for families of cycles
 modulo algebraic equivalence.

 Let $\pi: \mathcal{Y}\rightarrow B$ be a smooth projective morphism, where
 $B$ is smooth quasi-projective, and
 $\mathcal{Z}\subset \mathcal{Y}$ be a cycle of codimension $n$.
 Following Nori \cite{nori}, consider the
 cycle class
 $$[\mathcal{Z}]^{n,n}\in H^n(\mathcal{Y},\Omega_{\mathcal{Y}}^n)$$
 with restriction
 at $t\in B$:
 $$[\mathcal{Z}]^{n,n}_{t}\in H^n(\mathcal{Y}_t,{\Omega_{\mathcal{Y}}^n}_{\mid Y_t}).$$
 We have the following lemma, which will be useful later on. It  is not completely immediate because
   we are not working on a quasiprojective variety.
 \begin{lemm}\label{lealphann}
 The class $[\mathcal{Z}]^{n,n}$ is determined by the Betti cohomology class
  $[\mathcal{Z}]\in H^{2n}(\mathcal{Y},\mathbb{Q})$.
  \end{lemm}
  {\bf Proof.} The class $[\mathcal{Z}]$ belongs to the pure weight $2n$ part $W_{pure}H^{2n}(\mathcal{Y},\mathbb{Q})$
   of
   the Deligne mixed Hodge structure on
   $H^{2n}(\mathcal{Y},\mathbb{Q})$, which is defined as the image
   of the restriction map
   $$H^{2n}(\overline{\mathcal{Y}},\mathbb{Q})\rightarrow
   H^{2n}({\mathcal{Y}},\mathbb{Q})$$
   for any smooth projective completion
   $\overline{\mathcal{Y}}$ of $\mathcal{Y}$.
It is a Hodge class in $W_{pure}H^{2n}(\mathcal{Y},\mathbb{Q})$,
which means that it comes from a Hodge class $\alpha$ on
$\overline{\mathcal{Y}}$. This Hodge class can now be seen via the
Hodge decomposition as a class $\alpha^{n,n}\in
H^n(\overline{\mathcal{Y}},\Omega_{\overline{\mathcal{Y}}}^n)$,
which restricts to a class in
$H^n({\mathcal{Y}},\Omega_{{\mathcal{Y}}}^n)$. If we change $\alpha$
to $\alpha'$, the difference $\alpha-\alpha'$ is a Hodge class in
$\text{Ker}\,(H^{2n}(\overline{\mathcal{Y}},\mathbb{Q})\rightarrow
H^{2n}(\mathcal{Y},\mathbb{Q}))$. It is known (cf. \cite[Proof of
Lemma 2.22]{voisinweyl} or \cite[Corollary 7.3]{PS}) that
 if a Hodge class $\beta$ of degree $2n$ on a smooth projective
  compactification $\overline{\mathcal{Y}}$ of $\mathcal{Y}$ vanishes on $\mathcal{Y}$, then
  there exists a Hodge class $\beta'$ of degree $2n-2$
  on a desingularization $\overline{D}$ of $D:=\overline{\mathcal{Y}}\setminus \mathcal{Y}$ such that $j_*\beta'=\beta$,
   where $j:\overline{D}\rightarrow \overline{\mathcal{Y}}$ is the natural morphism.
  It follows that
  $j_*{\beta'}^{n-1,n-1}=\beta^{n,n}$, with
  ${\beta'}^{n-1,n-1}\in H^{n-1}(\overline{D},\Omega_{\overline{D}}^{n-1})$ and thus $\beta^{n,n}$ vanishes in
  $H^n(\mathcal{Y},\Omega_{\mathcal{Y}}^n)$.
  Applying this to $\beta=\alpha-\alpha'$, we conclude that $\alpha^{n,n}={\alpha'}^{n,n}$ in $H^n(\mathcal{Y},\Omega_{\mathcal{Y}}^n)$.
  \cqfd

 The class $[\mathcal{Z}]^{n,n}$  was again  considered and studied in \cite{voisinvar}, in the case of
 a family of $0$-cycles in a family of surfaces. Its usefulness for the study
 of the algebraic cycle $\mathcal{Z}$ comes from the following easy fact
 (see \cite[Theorem 10.19]{voisinbook}) :
 \begin{lemm} \label{lemmarat} If for any $t\in B$, the restricted cycle
 $\mathcal{Z}_t:=\mathcal{Z}_{\mid \mathcal{Y}_t}$ is rationally equivalent to
 $0$, then there is a dense Zariski open set $U\subset B$ such that $\mathcal{Z}_U$ is of torsion
 in $CH^n(\mathcal{Y}_U)$, hence in particular
  $$[\mathcal{Z}]^{n,n}=0\,\,{\rm in}\,\, H^n(\mathcal{Y}_U,{\Omega_{\mathcal{Y}_U}^n}),\,\,\,\,\,\,
  [\mathcal{Z}]^{n,n}_{t}=0\,\,{\rm in}\,\, H^n(\mathcal{Y}_t,{\Omega_{\mathcal{Y}}^n}_{\mid \mathcal{Y}_t}),\,t\in U.$$
  \end{lemm}
  In \cite{voisinvar}, \cite{voisinicm} (see also \cite[5.2.2]{voisinbook}), it was furthermore explained how to
  study the non-vanishing of the class $[\mathcal{Z}]^{n,n}$  by considering the ``Leray
  filtration''
  on ${\Omega_{\mathcal{Y}}^n}$:
  \begin{eqnarray}
  \label{leray}
  L^s{\Omega_{\mathcal{Y}}^n}:=\pi^*\Omega_{B}^s\wedge {\Omega_{\mathcal{Y}}^{n-s}}.
  \end{eqnarray}
  The associated graded is
  $$Gr_L^s{\Omega_{\mathcal{Y}}^n}=\pi^*\Omega_{B}^s\otimes
  {\Omega_{\mathcal{Y}/B}^{n-s}}.$$
  It thus follows that there is a
  spectral sequence abutting to  $R^n\pi_*{\Omega_{\mathcal{Y}}^n}$ :
  \begin{eqnarray}\label{spectral}
  E_1^{s,q}(\Omega_{\mathcal{Y}}^n,L)=\Omega_{B}^s\otimes R^{s+q}\pi_* {\Omega_{\mathcal{Y}/B}^{n-s}}\Rightarrow R^{s+q}\pi_*{\Omega_{\mathcal{Y}}^n}.
  \end{eqnarray}
  Of course, if $B$ is affine, we have
$$H^n(\mathcal{Y},\Omega_{\mathcal{Y}}^{n})=H^0(B,R^n\pi_*\Omega_{\mathcal{Y}}^n)$$
and thus the vanishing of $[\mathcal{Z}]^{n,n}$ is equivalent to the
vanishing of $[\mathcal{Z}]^{n,n}_t$ for any $t\in B$.

  One verifies (see \cite[5.1]{voisinbook}), as a consequence of the Katz-Oda
  description of the Gauss-Manin connection,
  that the first differential $d_1$ are  the Griffiths $\overline{\nabla}$-maps induced by
  Gauss-Manin connection and Griffiths transversality :
  \begin{eqnarray}\label{eqrel1}
  d_1=\overline{\nabla}_{s,n-s,s+q}:\Omega_{B}^s\otimes \mathcal{H}^{n-s,s+q}\rightarrow \Omega_{B}^{s+1}\otimes \mathcal{H}^{n-s-1,s+q+1},
  \end{eqnarray}
  where the  bundles
  $\mathcal{H}^{p,q}=R^q\pi_*\Omega_{\mathcal{Y}/B}^p$ are the Hodge
  bundles with fiber at $t\in B$ the spaces
  $H^q(\mathcal{Y}_t,\Omega_{\mathcal{Y}_t}^p)$.

  This spectral sequence degenerates at $E_2$ (see \cite{ikeda}, \cite{espa}), as it is the case
  for the topological Leray spectral sequence, and by the same  argument  as in   \cite{deligne}.

  Assume now  that our  cycle $\mathcal{Z}$ has its class in
  $L^sH^n(\mathcal{Y},\Omega_{\mathcal{Y}}^{n})$,
  which is equivalent, if $B$ is affine, to the fact that
  $[\mathcal{Z}]^{n,n}_{t}$ belongs to  the subspace $ L^sH^n(Y_t,{\Omega_{\mathcal{Y}}^n}_{\mid Y_t})$
  for any $t\in B$. Then it has an
  ``infinitesimal invariant'' (see \cite{voisinicm}, \cite[5.2.2]{voisinbook})
  \begin{eqnarray}\label{inf}\delta[\mathcal{Z}]\in E_\infty^{s,n-s}(\Omega_{\mathcal{Y}}^n,L)=E_2^{s,n-s}(\Omega_{\mathcal{Y}}^n,L)
  \\
  \nonumber
  =
 \text{ Ker}\,\overline{\nabla}_{s,n-s,n}/\text{Im}\,\overline{\nabla}_{s-1,n-s+1,n-1}
  \end{eqnarray}
  where
  $$\overline{\nabla}_{s,n-s,n}:\Omega_{B}^s\otimes \mathcal{H}^{n-s,n}\rightarrow \Omega_{B}^{s+1}\otimes \mathcal{H}^{n-s-1,n+1}$$
  and
  $$\overline{\nabla}_{s-1,n-s+1,n-1}:\Omega_{B}^{s-1}\otimes \mathcal{H}^{n-s+1,n-1}\rightarrow \Omega_{B}^{s}\otimes \mathcal{H}^{n-s,n},$$
are the $\overline{\nabla}$-maps of (\ref{eqrel1}).

Let $U\subset B$ be the dense Zariski open set where the coherent
sheaf $E_\infty^{s,n-s,n}$ is locally free.
  The non vanishing of the restricted infinitesimal invariant
  $$\delta[\mathcal{Z}]_t\in\text{ Ker}\,\overline{\nabla}_{s,n-s,n,t}/\text{Im}\,\overline{\nabla}_{s-1,n-s+1,n-1,t}$$
  at some point $t\in U$ implies that
  the class $[\mathcal{Z}]^{n,n}$ is nonzero on any dense
  Zariski open subset of $U$, and thus by Lemma \ref{lemmarat}, the cycle $\mathcal{Z}_t$ is not rationally equivalent to
  $0$ at a very general point $t'\in B$.

Let us comment on the assumption made on $\mathcal{Z}^{n,n}$, namely
that it belongs to $L^sH^n(\mathcal{Y},\Omega_{\mathcal{Y}}^{n})$,
and relate it to the assumption made in Theorem \ref{theoinfiintro}
that $ [\mathcal{Z}]\in L^sH^{2n}(\mathcal{Y},\mathbb{Q})$.

 \begin{prop} \label{proleray2avril} Let $\mathcal{Z}$ be a codimension $n$ cycle on $\mathcal{Y}$ such
 that the Betti cohomology class $[\mathcal{Z}]$ belongs to
 $L^sH^{2n}(\mathcal{Y},\mathbb{Q})$. Then, shrinking  the base $B$ if necessary,
 the Dolbeault cycle class $[\mathcal{Z}]^{n,n}$ belongs to
 $L^sH^n(\mathcal{Y},\Omega_\mathcal{Y}^n)$.
 \end{prop}
 {\bf Proof.} We use the fact (which will be proved later on, cf. the proof of Lemma \ref{lemmainter})
 that there exist Hodge classes $\delta_i$ of degree
 $2d$, $d:=\text{dim}\,\mathcal{Y}/B$, on $\mathcal{Y}\times_B\mathcal{Y}$,
 (where by Hodge classes, we mean Hodge classes in the pure part $W_{pure}H^{2d}(\mathcal{Y}\times_B\mathcal{Y},\mathbb{Q})$) satisfying
 the following properties:
 The classes $\delta_i$ act as K\"{u}nneth projectors, which means the following :
 The induced morphisms $\delta_{i*}:R\pi_*\mathbb{Q}\rightarrow R\pi_*\mathbb{Q}$
 factors as
 $$R\pi_*\mathbb{Q}\rightarrow R^i\pi_*\mathbb{Q}[-i]\rightarrow R\pi_*\mathbb{Q},$$
 where both  maps induce the identity on
 cohomology of degree $i$ (and $0$ in the other degrees). Furthermore, the sum
 $\sum_i\delta_i$ is equal to the class of the relative diagonal.

 When the fibres $\mathcal{Y}_t$ admit a Chow-K\"{u}nneth decomposition (cf. \cite{murre}), with components
 $\Delta_{i,t}$,
 which is the case of an abelian fibration,
 we can construct the $\delta_i$ over a generically finite cover $B'$ of $B$,
 by spreading out the $\Delta_{i,t}$ (which needs to make a base change) and taking
 the cohomology class $\delta_i$ of the spread-out cycles $\Delta_i$.

 We now let  the $\delta_i$ act on the classes $[\mathcal{Z}]$ and $[\mathcal{Z}]^{n,n}$.
 For the second action, we use the fact that $\delta_i$ is a Hodge class on $\mathcal{Y}\times_B\mathcal{Y}$
 hence has a natural image in $H^d(\mathcal{Y}\times_B\mathcal{Y},\Omega_{\mathcal{Y}\times_B\mathcal{Y}}^d)$.
 As we have $\sum_i\delta_i=[\Delta_{\mathcal{Y}/B}]$,
 we get
 \begin{eqnarray}
 \label{eqn2avril}\sum_i\delta_{i*}=Id: H^{2n}(\mathcal{Y},\mathbb{Q})\rightarrow
 H^{2n}(\mathcal{Y},\mathbb{Q}).
\end{eqnarray}
 Note also that the actions $\delta_{i*}$ are compatible with the Leray filtration on
 $H^{2n}(\mathcal{Y},\mathbb{Q})$ and the filtration $L$ on $H^{n}(\mathcal{Y},\Omega_{\mathcal{Y}}^n)$.
 Looking at the way the $\delta_i$ operate on the spectral sequence associated to the
 filtration $L$ on $\Omega_{\mathcal{Y}}^n$, one finds that
 $$\delta_{i*}[\mathcal{Z}]^{n,n}\in L^{2n-i}H^n(\mathcal{Y},\Omega_{\mathcal{Y}}^n).$$

 We write now using (\ref{eqn2avril})
 $$[\mathcal{Z}]=\sum_i\delta_{i*}[\mathcal{Z}].$$

 By Lemma \ref{lealphann}, we get that
if the Hodge class $\delta_{i*}[\mathcal{Z}]$ vanishes in the direct
 $H^{2n}(\mathcal{Y},\mathbb{Q})$,
 then its Dolbeault counterpart $\delta_{i*}[\mathcal{Z}]^{n,n}$ vanishes
 in $H^n(\mathcal{Y},\Omega_{\mathcal{Y}}^n)$.

 The assumption that $[\mathcal{Z}]\in
 L^sH^{2n}(\mathcal{Y},\mathbb{Q})$ is equivalent to the fact that
$\delta_{i*}[\mathcal{Z}]=0$ for $i>2n-s$, because
$\text{Im}\,\delta_{i*}:H^{2n}(\mathcal{Y},\mathbb{Q})\rightarrow
H^{2n}(\mathcal{Y},\mathbb{Q})$ is a direct summand contained in
$L^{2n-i}H^{2n}(\mathcal{Y},\mathbb{Q})$ and isomorphic to
$Gr_L^{2n-i}H^{2n}(\mathcal{Y},\mathbb{Q})$. But then
$\delta_{i*}[\mathcal{Z}]^{n,n}$ vanishes for $i>2n-s$, that is,
$$[\mathcal{Z}]^{n,n}=\sum_i\delta_{i*}[\mathcal{Z}]^{n,n}=\sum_{i\leq 2n-s}\delta_{i*}[\mathcal{Z}]^{n,n}\in
 \sum_{s\leq 2n-i} L^{2n-i}H^n(\mathcal{Y},\Omega_{\mathcal{Y}}^n)=
 L^{s}H^n(\mathcal{Y},\Omega_{\mathcal{Y}}^n).$$

 \cqfd
 The name ``infinitesimal invariant'' comes from the fact, proved in \cite{voisincras}, that
 for $s=1$, the infinitesimal invariant is the Griffiths infinitesimal
 invariant  (cf. \cite{griffithsIVHS})
 obtained by differentiating the normal function associated to
 $\mathcal{Z}$.
 \subsection{Infinitesimal invariants for cycles modulo algebraic equivalence \label{secalgeqinf}}
Our main result in this section is the construction of an
infinitesimal invariant for cycles {\it modulo algebraic
equivalence} and the proof of Theorem \ref{theoinfiintro}. Let $\pi:
\mathcal{Y}\rightarrow B$ be a smooth projective morphism and
$\mathcal{Z}\subset \mathcal{Y}$ be a cycle of codimension $n$. Note
that Lemma \ref{lemmarat} is not true for cycles modulo algebraic
equivalence. However, we have the following Proposition
\ref{theoinfi} which will allow us to extract from the
$\delta[\mathcal{Z}]$ defined in previous section an infinitesimal
invariant associated to cycles modulo algebraic equivalence.

We will use the fact that the maps
$$\overline{\nabla}_{s,n-s,n}:\Omega_{B}^{s}\otimes
\mathcal{H}^{n-s,n}\rightarrow \Omega_{B}^{s+1}\otimes
\mathcal{H}^{n-s-1,n+1}$$  satisfy the relation:
\begin{eqnarray}
\label{eqnew5avril}\overline{\nabla}_{s,n-s,n}(\alpha\wedge\sigma)=-\alpha\wedge\overline{\nabla}_{s-1,n-s,n}(\sigma)
\end{eqnarray}
for $\alpha$ a section of $\Omega_{B}$, $\sigma$ a section of
$\Omega_{B}^{s-1}\otimes \mathcal{H}^{n-s,n}$. This is proved as
follows :  for $\alpha$, $\sigma$ as above, let $\tilde{\sigma}$ be
a section of $F^{n-s}\mathcal{H}^{2n-s}$ which lifts $\sigma$. Then
by definition,
$$\overline{\nabla}_{s,n-s,n}(\alpha\wedge
\sigma)=\nabla(\alpha\wedge\tilde{\sigma})\,\,\text{mod.}\,\,\Omega_{B}^{s+1}\otimes
F^{n}\mathcal{H}^{2n-s}.$$
 By  Leibniz rule, this is also equal to
 $-\alpha\wedge \nabla\tilde{\sigma}+d\alpha\wedge \tilde{\sigma}\,\,\text{mod.}\,\,\Omega_{B}^{s+1}\otimes
F^{n}\mathcal{H}^{2n-s}$. As $d\alpha\wedge \tilde{\sigma}$ is a
section of $\Omega_{B}^{s+1}\otimes F^{n}\mathcal{H}^{2n-s}$, we get
$$\overline{\nabla}_{s,n-s,n}(\alpha\wedge
\sigma)= -\alpha\wedge
\nabla\tilde{\sigma}\,\,\text{mod.}\,\,\Omega_{B}^{s+1}\otimes
F^{n}\mathcal{H}^{2n-s}.$$ It follows from (\ref{eqnew5avril}) that
\begin{eqnarray} \label{wedge}\Omega_{B}\wedge
\text{Ker}\,\overline{\nabla}_{s-1,n-s,n}\subset
\text{Ker}\,\overline{\nabla}_{s,n-s,n},
 \end{eqnarray}
 as implicitly stated in  Theorem \ref{theoinfiintro}. We have now:
\begin{prop} \label{theoinfi}  Assume that for any $t\in B$,
(a multiple of) $\mathcal{Z}_t$ is algebraically equivalent to $0$
and furthermore $[\mathcal{Z}]^{n,n}\in
L^sH^n(\mathcal{Y},{\Omega_{\mathcal{Y}}^n})$. Then there is a
Zariski open set $U\subset B$ such that over $U$, the infinitesimal
invariant $\delta[\mathcal{Z}]^{n,n}$ belongs to the image of
$\Omega_{B}\wedge \text{Ker}\,\overline{\nabla}_{s-1,n-s,n}$ in
$$
\frac{\text{Ker}\,\overline{\nabla}_{s,n-s,n}}{\text{Im}\,\overline{\nabla}_{s-1,n-s+1,n-1}}.$$

\end{prop}
The following lemma is useful as it explains the simple form that
 the space of infinitesimal invariants in Theorem
\ref{theoinfiintro} takes.
\begin{lemm}  \label{leinclusion}  For $s\geq 2$, the
space $\Omega_{B}\wedge
\text{Ker}\,\overline{\nabla}_{s-1,n-s,n}\subset
\text{Ker}\,\overline{\nabla}_{s,n-s,n}$ contains
$\text{Im}\,\overline{\nabla}_{s-1,n-s+1,n-1}$.

It follows that the quotient
$$\frac{\text{Ker}\,\overline{\nabla}_{s-1,n-s,n}}{\text{Im}\,\overline{\nabla}_{s-1,n-s+1,n-1}+\Omega_{B}\wedge
\text{Ker}\,\overline{\nabla}_{s-1,n-s,n}}$$ is equal to
$$\frac{\text{Ker}\,\overline{\nabla}_{s-1,n-s,n}}{\Omega_{B}\wedge
\text{Ker}\,\overline{\nabla}_{s-1,n-s,n}}.$$
\end{lemm}
{\bf Proof.} Indeed, for $s\geq2$, we have $s-1\geq 1$ and thus for
a
 section $\eta$ of $ \Omega_{B}^{s-1}\otimes
 \mathcal{H}^{n-s+1,n-1}$,
 we can write locally $\eta=\sum_i\alpha_i\wedge\eta_i$, where
 $\alpha_i$ are sections of $\Omega_B$ and $\eta_i$ are sections of
 $ \Omega_{B}^{s-2}\otimes
 \mathcal{H}^{n-s+1,n-1}$.
 By (\ref{eqnew5avril}),
 $$\overline{\nabla}_{s-1,n-s+1,n-1}(\eta)=-\sum_i\alpha_i\wedge
 \overline{\nabla}_{s-2,n-s+1,n-1}(\eta_i).$$
 This proves the result since
 $\overline{\nabla}_{s-2,n-s+1,n-1}(\eta_i)\in
 \text{Ker}\,\overline{\nabla}_{s-1,n-s,n}$.

 \cqfd

The proof of Proposition \ref{theoinfi}  is postponed to the end of
the section. We first show how it implies Theorem
\ref{theoinfiintro}.

\vspace{0.5cm}

 {\bf Proof of
Theorem \ref{theoinfiintro}.} Starting from a cycle $\mathcal{Z}$
with class $[\mathcal{Z}]\in L^sH^{2n}(\mathcal{Y},\mathbb{Q})$, we
know by Proposition \ref{proleray2avril} that up to shrinking $B$ if
necessary, the class $[\mathcal{Z}]^{n,n}$ belongs to
$L^sH^n(\mathcal{Y},{\Omega_{\mathcal{Y}}^n})$. We thus have the
infinitesimal invariant $\delta[\mathcal{Z}]$ and its restrictions
$\delta[\mathcal{Z}]_{t}\in
Gr_L^sH^n(\mathcal{Y}_t,\Omega_{\mathcal{Y}\mid\mathcal{Y}_t}^n)$
for any point $t\in B$. Let $U$ be the dense Zariski open subset of
$B$ where the sheaf
$\frac{\text{Ker}\,\overline{\nabla}_{s,n-s,n}}{\Omega_{B,b}\wedge
\text{Ker}\,\overline{\nabla}_{s-1,n-s,n}}$ is locally free. Then
the non vanishing of the infinitesimal invariant
$\delta[\mathcal{Z}]_{t}\in
\text{Ker}\,\overline{\nabla}_{s,n-s,n,t}/\text{Im}\,\overline{\nabla}_{s-1,n-s-1,n-1,t}$
modulo the image of $\text{Im}\,\Omega_{B,t}\otimes
\text{Ker}\,\overline{\nabla}_{s-1,n-s,n,t}$ by the wedge product
map implies that for any dense Zariski open subset $U'\subset U$,
$\delta[\mathcal{Z}]_{\mid U'}$ does not vanish modulo
$\text{Im}(\Omega_{U'}\wedge
\text{Ker}\,\overline{\nabla}_{s-1,n-s,n})$,
 hence it follows from  Proposition \ref{theoinfi} that the
cycle $\mathcal{Z}_t$ is not algebraically equivalent to $0$ for
very general $t\in U$.

Denoting by $\delta[Z]_{alg}$, resp. $\delta[\mathcal{Z}]_{t,alg}$
the
 image of $\delta[\mathcal{Z}]$ in
 $\text{Ker}\,\overline{\nabla}_{s-1,n-s,n}/\Omega_U\wedge
 \text{Ker}\,\overline{\nabla}_{s-1,n-s,n}$,
 (resp.  $\delta[\mathcal{Z}]_t$ in $\text{Ker}\,\overline{\nabla}_{s-1,n-s,n,t}/\Omega_{B,t}\wedge
 \text{Ker}\,\overline{\nabla}_{s-1,n-s,n,t}$), we constructed the
 desired
  infinitesimal invariant for cycles modulo algebraic equivalence
  satisfying the conclusion of Theorem \ref{theoinfiintro} and Corollary
  \ref{coroinfiintro}.

\cqfd

Let us first make a few remarks.
\begin{rema} {\rm In \cite{nori}, Nori uses the cycle class $[\mathcal{Z}]^{n,n}$ to conclude
that a certain cycle is not algebraically equivalent to $0$ at a
general point $t\in B$, assuming its cohomology class $\mathcal{Z}$
does not vanish on the total space. However, he does not introduce
the infinitesimal invariant above (although the infinitesimal
computation is hidden in his arguments), for the following reason.
In his situation, the kernel of the map
$$\overline{\nabla}_{s-1,n-s,n}: \Omega_{B}^{s-1}\otimes \mathcal{H}^{n-s,n}\rightarrow \Omega_{B}^{s}\otimes \mathcal{H}^{n-s-1,n+1}$$
 is equal to the image of
the map
$$ \overline{\nabla}_{s-2,n-s+1,n-1}: \Omega_{B}^{s-2}\otimes \mathcal{H}^{n-s+1,n-1}\rightarrow \Omega_{B}^{s}\otimes \mathcal{H}^{n-s,n},$$
and it thus follows that
$$\Omega_{B}\wedge \text{Ker}\,\overline{\nabla}_{s-1,n-s,n}\subset \text{Im}\,(\overline{\nabla}_{s-1,n-s+1,n-1}: \Omega_{B}^{s-1}\otimes
 \mathcal{H}^{n-s+1,n-1}\rightarrow \Omega_{B}^{s}\otimes \mathcal{H}^{n-s,n}).$$
This  is why  Nori can use the previously defined invariant
$\delta[\mathcal{Z}]$. The same remark applies to the paper of
Fakhruddin who uses similar methods applied to generic abelian
varieties. }
\end{rema}

Our second  remark concerns the consistency of the existence of the
above defined invariant and the conjecture made in
\cite{voisinpura}. We first recall what is this conjecture: Here we
can consider for the Bloch-Beilinson filtration on Chow groups with
rational coefficients
 any of the filtrations mentioned in the introduction, for example
the one defined by Shuji Saito in \cite{ssaito}, or, in the case of
an abelian variety, the one which is induced by the Beauville
decomposition as in (\ref{eqbeaudec}). Let us denote by $F$ such a
filtration. These filtrations are conjectured to satisfy the crucial
property:
$$F^{k+1}\text{CH}^k(Y)_\mathbb{Q}=0.$$
These filtrations also induce  similar filtrations on
Chow groups modulo algebraic equivalence.

We made in \cite{voisinpura} (and explored the consequences of)  the following conjecture
(which generalizes Nori's conjecture  \cite{nori} on $CH^2$):
\begin{conj} We have $$F^{k}\text{CH}^k(Y)/ alg=0.$$
\end{conj}

We want to observe now the following fact, which  is consistent with the above conjecture
\begin{lemm} Let $\mathcal{Z}\subset \mathcal{Y}\rightarrow B$ be a  cycle
of codimension $n$, such that $[\mathcal{Z}]\in
L^nH^n(\mathcal{Y},\Omega_{\mathcal{Y}}^n)$. Then for any $t\in B$,
the infinitesimal invariant
$$\delta[\mathcal{Z}]_{t,alg}\in \frac{\text{Ker}\,\overline{\nabla}_{n,0,n,t}}{\Omega_{B,t}\wedge
\text{Ker}\,\overline{\nabla}_{n-1,0,n,t}}
$$
vanishes.
\end{lemm}
{\bf Proof.} This is obvious because the group we are looking at is
the following (assuming $n\geq 2$, see Lemma \ref{leinclusion}):
$$\frac{Ker\,(\overline{\nabla}_{n,0,n}: \Omega_{B,t}^n\otimes H^{0,n}(\mathcal{Y}_t)\rightarrow
\Omega_{B_t}^{n+1}\otimes
H^{-1,n+1}(\mathcal{Y}_t))}{\Omega_{B,t}\wedge
Ker\,(\overline{\nabla}_{n-1}: \Omega_{B,t}^{n-1}\otimes
H^{0,n}(\mathcal{Y}_t)\rightarrow \Omega_{B,t}^n\otimes
H^{-1,n+1}(\mathcal{Y}_t))}.$$ Hence, because
$H^{-1,n+1}(\mathcal{Y}_t)=0$, the $\overline{\nabla}$-maps  vanish
in this range, and our group is in fact
$$\frac{\Omega_{B,t}^n\otimes H^{0,n}(\mathcal{Y}_t)}{\Omega_{B,t}\wedge \Omega_{B,t}^{n-1}\otimes H^{0,n}(\mathcal{Y}_t)}$$
which is $0$ because $n>0$. \cqfd
  Our last remark concerns the assumption
that $[\mathcal{Z}]\in L^sH^{2n}(\mathcal{Y},\mathbb{Q})$ and will
be useful for our main application. There are two cases where the
Leray level of an algebraic cycle in a family of varieties is easy
to compute:

1) The family $\mathcal{Y}$ is a family of smooth
 ample  hypersurfaces or complete intersections
in projective space (or more generally, any ambient variety $X$
whose rational cohomology is made of classes of algebraic cycles).
Then the hard Lefschetz theorem in this case says that the
interesting part of the cohomology of $\mathcal{Y}_t$ is supported
in degree $m:=dim\,\mathcal{Y}_t$. Given any cycle $\mathcal{Z}$ of
codimension $n<\text{dim}\,\mathcal{Y}_t/2$ in $ \mathcal{Y}$, we
can write (at least with rational coefficients)
$$\mathcal{Z}=\mathcal{Z}'+\Gamma$$
where $\Gamma$ is the restriction of a cycle coming from $X$ and
$\mathcal{Z}'_{\mid Y_t}$ is cohomologous to $0$. Looking more
precisely at the shape of the Leray spectral sequence of
$\mathcal{Y}\rightarrow B$ and applying the Lefschetz hyperplane
section theorem, one finds that on an adequate dense Zariski open
set $U\subset B$, the Dolbeault cohomology class of
$\mathcal{Z}'_{\mid \mathcal{Y}_U}$ belongs to
$L^{2n-d}H^{2n}(\mathcal{Y}_U,\mathbb{Q})$, where
$d=\text{dim}\,\mathcal{Y}_t$.

\vspace{0,5cm}

2)  The second case where the Leray level is easy to compute is the
case of a cycle in   a family of abelian varieties
$\mathcal{A}\rightarrow B$. In this case we have the following fact:
\begin{lemm} \label{lemmespecab} If the cycle $\mathcal{Z}$ satisfies
the condition that
$$Z_t\in \text{CH}^n(\mathcal{A}_t)_s$$
for any $t\in B$, then there is a dense Zariski open set $U\subset
B$ such that $[\mathcal{Z}]_{U}\in
L^sH^{2n}(\mathcal{A}_U,\mathbb{Q})$.

\end{lemm}
{\bf Proof.} Indeed, under our assumptions, Lemma \ref{lemmarat},
and the relative Beauville decomposition (see \cite{demu}) show that
there exists a dense Zariski open set $U\in B$ such that for any
$k\in \mathbb{Z}^*$,
\begin{eqnarray}\label{1}\mu_k^*\mathcal{Z}_U= k^{2n-s}\mathcal{Z}\,\,{\rm in} \,\,
\text{CH}^n(\mathcal{A}_U).
\end{eqnarray}
 It thus follows that
the class $[\mathcal{Z}]\in\,\, H^{2n}(\mathcal{A}_U,\mathbb{Q})$
satisfies for any $k$:
\begin{eqnarray}\label{2}\mu_k^*[\mathcal{Z}]= k^{2n-s}[\mathcal{Z}] \,\,{\rm in}\,\, H^{2n}(\mathcal{A}_U,\mathbb{Q}).
\end{eqnarray}
But there is a canonical Deligne decomposition
\begin{eqnarray}
\label{eqdedel}R\pi_*\mathbb{Q}=\oplus R^i\pi_*\mathbb{Q}[-i]
\end{eqnarray} such that for any $i$,
 $$\mu_k^*\circ\pi_i=k^i\pi_i :R\pi_*\mathbb{Q}\rightarrow R\pi_*\mathbb{Q},$$
 where the $\pi_i$ are the projectors associated to the decomposition
 (\ref{eqdedel}) (cf. \cite{demu}). It follows that
\begin{eqnarray}\label{3} L^sH^{2n}(\mathcal{A}_U,\mathbb{Q})=\oplus_{i\geq s}
 H^{2n}(\mathcal{A}_U,\mathbb{Q})_{i},
 \end{eqnarray}
 where $$H^{2n}(\mathcal{A}_U,\mathbb{Q})_{i}
 =\{\alpha\in H^{2n}(\mathcal{A}_U,\mathbb{Q}),\,\mu_k^*[\mathcal{Z}]= k^{2n-i}[\mathcal{Z}] \,\,{\rm
 in}  \,\,
 H^{2n}(\mathcal{A}_U,\mathbb{Q})\,\,
 \text{for any}\,k\}.$$
 By (\ref{2}) and (\ref{3}), $[\mathcal{Z}]_{U}\in
 L^sH^{2n}(\mathcal{A}_U,\mathbb{Q})$.
  \cqfd

{\bf Proof of Proposition \ref{theoinfi}.} We first prove the result
in the case which will be useful for us, namely when
$\pi:\mathcal{Y}\rightarrow B$ is a family of abelian varieties. The
proof is easier in this case because, as the proof of lemma
\ref{lemmespecab} shows, not only the spectral sequence abutting to
$R^n\pi_*(\mathcal{Y},{\Omega_{\mathcal{Y}}^n})$ degenerates at
$E_2$ but there is a natural decomposition
\begin{eqnarray}\label{4}R^n\pi_*{\Omega_{\mathcal{Y}}^n}=\oplus_s
(R^n\pi_*{\Omega_{\mathcal{Y}}^n})_s,
 \end{eqnarray}
 where
$$ (R^n\pi_*{\Omega_{\mathcal{Y}}^n})_s=\{\alpha\in R^n\pi_*{\Omega_{\mathcal{Y}}^n},\,\mu_k^*\alpha=k^{2n-s}\alpha,\,\forall k\in\mathbb{Z}\}$$
and the decomposition (\ref{4}) is a splitting of the filtration
induced by the filtration $L$ on $\Omega_{\mathcal{Y}}$. We will
explain at the end of the proof how to modify the argument to make
it work in general.

So assume $\mathcal{Y}\rightarrow B$ is a family of abelian
varieties over a smooth affine variety $B$, and consider a
codimension $n$ cycle $\mathcal{Z}$ on $ \mathcal{Y}$. Assume that
for any $t\in B$, $\mathcal{Z}_t$ is algebraically equivalent to
$0$, and furthermore $[\mathcal{Z}]\in
L^sH^{2n}(\mathcal{Y},\mathbb{Q})$. Using the relative Beauville
decomposition of Deninger-Murre \cite{demu}, we can replace
$\mathcal{Z}$ by $\mathcal{Z}_s$, which does not change the class
$\delta[\mathcal{Z}]\in Gr_s^LH^n(\mathcal{Y},\Omega_{
\mathcal{Y}}^n)$ and it still satisfies the assumption that its
restrictions to the fibers is algebraically equivalent to $0$. So we
 assume from now on that $\mathcal{Z}=\mathcal{Z}_s$.

The definition of algebraic equivalence together with an elementary
argument involving countability of Chow varieties and uncountability
of $\mathbb{C}$ implies that over a Zariski dense open subset $U$ of
$ B$, there exist
\begin{enumerate}
\item A family of smooth non necessarily connected curves $\rho:\mathcal{C}\rightarrow
U$;
\item A divisor $\mathcal{D}$ of $\mathcal{C}$ which is homologous to $0$ on the fibers of
$\rho$;
\item A codimension $n$ correspondence  $\Gamma\subset \mathcal{C}\times_U\mathcal{Y}_U$ with
$\mathbb{Q}$-coefficients,

\end{enumerate}
such that
\begin{eqnarray} \label{eqcycles}\mathcal{Z}=\Gamma_*(\mathcal{D})\,\,{\rm in}\,\,\text{CH}^n(\mathcal{Y}_U).
\end{eqnarray}

Here we recall that, denoting by $p_1,\,p_2$ the two proper smooth projections from
$\mathcal{C}\times_U\mathcal{Y}_U$ to $\mathcal{C}$ and $\mathcal{Y}_U$ respectively,
$$\Gamma_*(\mathcal{D}):=p_{2*}(p_1^*\mathcal{D}\cdot \Gamma).$$

 The
equality of cycles (\ref{eqcycles}) provides now, using the natural
functoriality properties of the Dolbeault cycle class, the
corresponding equality of Dolbeault cohomology classes:
\begin{eqnarray} \label{eqclasses}[\mathcal{Z}_U]^{n,n}=p_{2*}([\Gamma]^{n,n}\cdot p_1^*[\mathcal{D}]^{1,1})\,\,{\rm in}\,\,
H^n(\mathcal{Y}_U,{\Omega_{\mathcal{Y}_U}^n}).
\end{eqnarray}
We use now again the relative Beauville decomposition for cycles in
the  family  $\mathcal{C}\times_U\mathcal{Y}\rightarrow
\mathcal{C}$. The cycle $\Gamma$ decomposes as
$$\Gamma=\sum_{p}\Gamma_p,$$
with $\mu_k^*\Gamma_p=k^{2n-p}\Gamma_p$ and accordingly, the class
$$[\Gamma]^{n,n}\in H^n(\mathcal{C}\times_U\mathcal{Y},{\Omega_{\mathcal{C}\times_U\mathcal{Y}}^n})$$
decomposes then into a direct sum
$$[\Gamma]^{n,n}=\oplus_p [\Gamma_p]^{n,n},$$
where $[\Gamma_p]^{n,n}\in
H^n(\mathcal{C}\times_U\mathcal{Y},{\Omega_{\mathcal{C}\times_U\mathcal{Y}}^n})$
satisfies the condition
\begin{eqnarray} \label{decclassesgamma}\mu_k^*[\Gamma_p]^{n,n}=k^{2n-p}[\Gamma_p]^{n,n},\,\forall k\in \mathbb{Z}^*.
\end{eqnarray}
Recall that $[\mathcal{Z}]^{n,n}_U$ belongs to
$H^n(\mathcal{Y}_U,{\Omega_{\mathcal{Y}_U}^n})_s$. Observing that
the $\mu_k^*$ maps are compatible with the pushforward map $p_{2*}$,
in the sense that

\begin{eqnarray} \label{decmu} p_{2*}\circ \mu_k^* =\mu_k^*\circ p_{2*},
\end{eqnarray}

we get  from (\ref{eqclasses}), (\ref{decclassesgamma}) and
(\ref{decmu})
 the following equality:
\begin{eqnarray} \label{eqclassestronq}[\mathcal{Z}_U]^{n,n}=p_{2*}([\Gamma_s]^{n,n}\cdot p_1^*[\mathcal{D}]^{1,1})\,\,{\rm in}\,\,
H^n(\mathcal{Y},{\Omega_{\mathcal{Y}}^n})_s\subset
L^sH^n(\mathcal{Y},{\Omega_{\mathcal{Y}}^n}),
\end{eqnarray}
and projecting everything modulo $L^{s+1}$, we get  :
\begin{eqnarray} \label{eqclassestronqpile}\delta[\mathcal{Z}]=p_{2*}(\delta[\Gamma_s]\cdot p_1^*\delta[\mathcal{D}])\,\,{\rm in}\,\,
Gr_L^sH^n(\mathcal{Y}_U,{\Omega_{\mathcal{Y}_U}^n})=E_2^{s,n-s}(\Omega_{\mathcal{Y}_U}^n,L).
\end{eqnarray}
Here $[\mathcal{D}]^{1,1}\in
H^1(\mathcal{C},{\Omega_{\mathcal{C}}})$, and as $\mathcal{D}$ is
cohomologous to $0$ on the fibers $\mathcal{C}_t$, $t\in U$, we have
$[\mathcal{D}]\in L^1 H^1(\mathcal{C},{\Omega_{\mathcal{C}}})$, so
\begin{eqnarray}
\label{eqdeltaD} \delta[\mathcal{D}]\in \Omega_{U}\otimes
\mathcal{H}^{0,1}_{\mathcal{C}}/\overline{\nabla}_{\mathcal{C}}\mathcal{H}^{1,0}_{\mathcal{C}},
\end{eqnarray}
where the bundles
$\mathcal{H}^{0,1}_{\mathcal{C}}=R^1\rho_*\mathcal{O}_\mathcal{C}$,
$\mathcal{H}^{1,0}_{\mathcal{C}}=R^0\rho_*\Omega_{\mathcal{C}/U}$
are the Hodge bundles of the family $\rho:\mathcal{C}\rightarrow U$
and $\nabla_\mathcal{C},\,\overline{\nabla}_\mathcal{C}$ denote the
corresponding Gauss-Manin connection and $\overline{\nabla}$-map.
(To avoid confusion, we will also write below
$\overline{\nabla}_\mathcal{Y}$, $\mathcal{H}^{p,q}_\mathcal{Y}$ for
the infinitesimal variation of Hodge structure of the family
$\mathcal{Y}$.)
 We will denote below by $\widetilde{\delta[\mathcal{D}]}$ any
lift of $\delta[\mathcal{D}]$ in $\Omega_{U}\otimes
\mathcal{H}^{0,1}_{\mathcal{C}}$.

 We now examine
$$\gamma:=\delta[\Gamma_s]\in
 H^{n}(\mathcal{C}\times_U\mathcal{Y},{\Omega_{\mathcal{C}\times_U\mathcal{Y}_U}^n})_s=H^0(U,(R^n\pi'_*\Omega_{\mathcal{C}\times_U\mathcal{Y}_U}^n)_s),$$
  where
  $\pi':=\pi\circ p_2:\mathcal{C}\times_U \mathcal{Y}_U\rightarrow
 U$, and
 as before
 $$(R^n\pi'_*\Omega_{\mathcal{C}\times_U\mathcal{Y}_U}^n)_s=\{\alpha\in
 R^n\pi'_*\Omega_{\mathcal{C}\times_U\mathcal{Y}_U}^n,\,\mu_k^*\alpha=k^{2n-s}\alpha\}.$$

 We recall  that the class $\gamma$ is a section of
$$(R^n\pi'_*\Omega_{\mathcal{C}\times_U\mathcal{Y}_U}^n)_s,$$
with $\pi'=\rho\circ p_1=\pi\circ
p_2:\mathcal{C}\times_U\mathcal{Y}_U\rightarrow U$.
 Considering the exact sequence

$$0\rightarrow R^1\rho_*(R^{n-1}p_{1*}\Omega_{\mathcal{C}\times_U\mathcal{Y}_U}^n)
\rightarrow R^n\pi'_*\Omega_{\mathcal{C}\times_U\mathcal{Y}_U}^n
\rightarrow
\rho_*(R^np_{1*}\Omega_{\mathcal{C}\times_U\mathcal{Y}_U}^n)\rightarrow
0,
$$
and taking the $s$-th direct summand (where $\mu_k^*$ acts by
multiplication by $k^{2n-s}$), we get:
$$0\rightarrow
R^1\rho_*(R^{n-1}p_{1*}\Omega_{\mathcal{C}\times_U\mathcal{Y}_U}^n)_s
\rightarrow (R^n\pi'_*\Omega_{\mathcal{C}\times_U\mathcal{Y}_U}^n)_s
\rightarrow
\rho_*(R^np_{1*}\Omega_{\mathcal{C}\times_U\mathcal{Y}_U}^{n})_s\rightarrow
0.
$$
As $\Omega_{\mathcal{C}\times_U\mathcal{Y}_U}^n$ projects naturally
to $p_1^*\Omega_{\mathcal{C}/U}\otimes
p_2^*\Omega_{\mathcal{Y}}^{n-1}$, we get a morphism
$$\rho_*(R^np_{1*}\Omega_{\mathcal{C}\times_U\mathcal{Y}_U}^n)_s\rightarrow \mathcal{H}^{1,0}_\mathcal{C}\otimes
 (R^n\pi_*\Omega_{\mathcal{Y}_U}^{n-1})_s.$$
 As was explained in section \ref{subsecrat},
$(R^n\pi_*\Omega_{\mathcal{Y}_U}^{n-1})_s$ is computed as
$$\frac{\text{Ker}\,(\overline{\nabla}_{\mathcal{Y},s-1,n-s,n}: \Omega_U^{s-1}\otimes
\mathcal{H}_{\mathcal{Y}}^{n-s,n}\rightarrow
\Omega_U^{s}\otimes\mathcal{H}_{\mathcal{Y}}^{n-s-1,n+1})}{\text{Im}\,(\overline{\nabla}_{\mathcal{Y},s-2,n-s+1,n-1}:
\Omega_U^{s-2}\otimes \mathcal{H}\mathcal{Y}^{n-s+1,n-1}\rightarrow
\Omega_U^{s-1}\otimes\mathcal{H}_\mathcal{Y}^{n-s,n})}.$$
 Combining these various
morphisms and identifications, we find that $\gamma$ provides a
class
$$\gamma_1\in\mathcal{H}^{1,0}_\mathcal{C}\otimes
\frac{\text{Ker}\,\overline{\nabla}_{\mathcal{Y},s-1,n-s,n}}{\text{Im}\,\overline{\nabla}_{\mathcal{Y},s-2,n-s+1,n-1}}.$$
Below, we will denote by $\tilde{\gamma}_1$ any lift of $\gamma_1$
in
$\mathcal{H}^{1,0}_\mathcal{C}\otimes\text{Ker}\,\overline{\nabla}_{\mathcal{Y},s-1,n-s,n}\subset
\mathcal{H}^{1,0}_\mathcal{C}\otimes\Omega_U^{s-1}\otimes\mathcal{H}_{\mathcal{Y}}^{n-s,n}$.

The proof of Proposition \ref{theoinfi} in the case of a family of
abelian varieties is then an immediate consequence of the following
Lemma \ref{lemmacontra}. Indeed, it says in particular that, under
the assumptions of Proposition \ref{theoinfi},
$\delta[\mathcal{Z}]\in
\frac{\text{Ker}\,\overline{\nabla}_{\mathcal{Y},s,n-s,n}}{\text{Im}\,\overline{\nabla}_{\mathcal{Y},s-1,n-s+1,n-1}}$
belongs to the image, under the wedge product map, of the sheaf
$\Omega_U\otimes\frac{
\text{Ker}\,\overline{\nabla}_{\mathcal{Y},s-1,n-s,n}}{\text{Im}\,\overline{\nabla}_{\mathcal{Y},s-2,n-s+1,n-1}}$.
\begin{lemm}\label{lemmacontra} With the same assumptions and notation as above, we have
\begin{eqnarray}
\label{formdesiree8avril}
\delta[\mathcal{Z}]=w(<\widetilde{\delta[\mathcal{D}]},\tilde{\gamma}_1>)\,\,\text{modulo}\,\,
{\text{Im}\,\overline{\nabla}_{\mathcal{Y},s-1,n-s+1,n-1}},
\end{eqnarray}
 where the
brackets mean that we use the duality  between
$\mathcal{H}^{0,1}_\mathcal{C}$ and $\mathcal{H}^{1,0}_\mathcal{C}$
to get a contraction map
$$<,>:(\Omega_U\otimes\mathcal{H}^{0,1}_\mathcal{C} )\otimes(\mathcal{H}^{1,0}_\mathcal{C}\otimes\Omega_U^{s-1}\otimes  \mathcal{H}_{\mathcal{Y}}^{n-s,n})
\rightarrow \Omega_U\otimes\Omega_U^{s-1}\otimes
\mathcal{H}_{\mathcal{Y}}^{n-s,n},$$
 and the map $w$ is induced by the wedge product map
$\Omega_U\otimes\Omega_U^{s-1}\rightarrow \Omega_U^s$.

\end{lemm}
Before giving the proof of (\ref{formdesiree8avril}), let us explain
why the projection of
$$w(<\widetilde{\delta[\mathcal{D}]},\tilde{\gamma}_1>)\in
\text{Ker}\,\overline{\nabla}_{\mathcal{Y},s,n-s,n}\subset\Omega_U^s\otimes
\mathcal{H}\mathcal{Y}^{n-s,n}$$ in the quotient
$\frac{\text{Ker}\,\overline{\nabla}_{\mathcal{Y},s,n-s,n}}{\text{Im}\,\overline{\nabla}_{\mathcal{Y},s-1,n-s+1,n-1}}$
does not depend on the choice of the lifts
$\widetilde{\delta[\mathcal{D}]},\,\tilde{\gamma}_1$. If we change
$\tilde{\gamma}_1$ by a $\overline{\nabla}_{\mathcal{Y}}$-exact
section of
$\Omega_U^{s-1}\otimes\mathcal{H}^{1,0}_\mathcal{C}\otimes\mathcal{H}_{\mathcal{Y}}^{n-s,n}$,
that is $(\tilde{\gamma}_1)'=
\tilde{\gamma}_1+\overline{\nabla}_{\mathcal{Y},s-2,n-s+1,n-1}(\gamma')$
for a section $\gamma'$ of
$\Omega_U^{s-2}\otimes\mathcal{H}^{1,0}_\mathcal{C}\otimes\mathcal{H}_{\mathcal{Y}}^{n-s+1,n-1}$,
then we clearly have:
$$w(<\delta[\widetilde{\mathcal{D}}],(\tilde{\gamma}_1)'>=w(<\widetilde{\delta[\mathcal{D}]},\tilde{\gamma}_1>)+
\overline{\nabla}_{\mathcal{Y},s-1,n-s+1,n-1}(w(<\delta[\widetilde{\mathcal{D}}],\gamma'>)),$$
and the difference is
$\overline{\nabla}_{\mathcal{Y},s-1,n-s+1,n-1}$-exact. On the other
hand, if we change $\widetilde{\delta[\mathcal{D}]}$ to
$$(\widetilde{\delta[\mathcal{D}]})'=\widetilde{\delta[\mathcal{D}]}+\overline{\nabla}_\mathcal{C}(\alpha)$$
where $\alpha$ is a section of $\mathcal{H}^{1,0}_\mathcal{C}$, then
we get
$$w(<(\widetilde{\delta[\mathcal{D}]})',\tilde{\gamma}_1>)=w(<\widetilde{\delta[\mathcal{D}]},\tilde{\gamma}_1>)+
w(<\overline{\nabla}_\mathcal{C}(\alpha),\tilde{\gamma}_1>)$$
$$=w(<\widetilde{\delta[\mathcal{D}]},\tilde{\gamma}_1>)-w(<\alpha,\overline{\nabla}_\mathcal{C}(\tilde{\gamma}_1)>).
$$
On the other hand, the class
$\overline{\nabla}_\mathcal{C}(\tilde{\gamma}_1)$, hence also the
contraction
$-w(<\alpha,\overline{\nabla}_\mathcal{C}(\tilde{\gamma}_1)>)$, is
$\overline{\nabla}_\mathcal{Y}$-exact for the following reason:
Recall that $\tilde{\gamma}_1$ is the projection of a
$\overline{\nabla}_{\mathcal{C}\times_B\mathcal{Y}}$-closed section
$\tilde{\gamma}$ of
$(R^n\pi'_*\Omega_{\mathcal{C}\times_B\mathcal{Y}/B}^n)_s$. This
section has a K\"{u}nneth decomposition
$\tilde{\gamma}=\tilde{\gamma}_1+\tilde{\gamma_2}+\ldots$ with
$\tilde{\gamma_2}\in
\Omega_U^{s-1}\otimes\mathcal{H}^{0,1}_\mathcal{C}\otimes
\mathcal{H}^{n-s+1,n-1}_{\mathcal{Y}}$. The condition
$\overline{\nabla}_{\mathcal{C}\times_U\mathcal{Y}_U}(\tilde{\gamma})=0$
implies then
$$\overline{\nabla}_\mathcal{C}(\tilde{\gamma}_1)=-\overline{\nabla}_\mathcal{Y}(\tilde{\gamma}_2).$$

 \vspace{0.5cm}

 {\bf Proof of
Lemma \ref{lemmacontra}.} We apply formula (\ref{eqclassestronq}).
We observe that the map
$$p_{2*}:R^{n+1}\pi'_*\Omega_{\mathcal{C}\times_U\mathcal{Y}_U}^{n+1}\rightarrow
R^{n}\pi_*\Omega_{\mathcal{Y}_U}^{n}$$
is obtained by composing the maps $I$ and $g$, where
$$g:R^{n+1}\pi'_*\Omega_{\mathcal{C}\times_U\mathcal{Y}_U}^{n+1}\rightarrow
R^{n+1}\pi'_*(p_1^*\Omega_{\mathcal{C}/U}\otimes
p_2^*\Omega_{\mathcal{Y}_U}^{n})$$ is  induced by  the morphism
$$\Omega_{\mathcal{C}\times_U\mathcal{Y}_U}^{n+1}\rightarrow
p_1^*\Omega_{\mathcal{C}/U}\otimes p_2^*\Omega_{\mathcal{Y}_U}^{n}$$
coming from the exact sequence:
$$ 0\rightarrow p_2^*\Omega_{\mathcal{Y}_U}
\rightarrow\Omega_{\mathcal{C}\times_U\mathcal{Y}_U}\rightarrow
p_1^*\Omega_{\mathcal{C}/U}\rightarrow 0$$

and
$$ I:R^{n+1}\pi'_*(p_1^*\Omega_{\mathcal{C}/U}\otimes p_2^*\Omega_{\mathcal{Y}_U}^{n})\rightarrow R^n\pi_*\Omega_{\mathcal{Y}_U}^{n}$$
is the integration over the fibre of $p_2$, also obtained
by writing $\pi'=\pi\circ p_2$ and sending
$R^{n+1}\pi'_*(p_1^*\Omega_{\mathcal{C}/U}\otimes p_2^*\Omega_{\mathcal{Y}_U}^{n})$ to
$$R^n\pi_*(R^1p_{2*}(p_1^*\Omega_{\mathcal{C}/U}\otimes p_2^*\Omega_{\mathcal{Y}_U}^{n}))=R^n\pi_*(\Omega_{\mathcal{Y}_U}^{n})$$
by the corresponding Leray spectral sequence.

We apply the morphism $g$ to the class
$p_1^*[\mathcal{D}]^{1,1}\cdot [\Gamma_s]^{n,n}$. As
$[\mathcal{D}]^{1,1}\in H^1(\mathcal{C},\Omega_\mathcal{C})$
vanishes in $ H^0(U,R^1\rho_*\Omega_{\mathcal{C}/U})$, we find that
\begin{eqnarray}
\label{eq8avrilprod} g(p_1^*[\mathcal{D}]^{1,1}\cdot
[\Gamma_s]^{n,n})=\widetilde{\delta[\mathcal{D}]}\cdot
g(\gamma)\,\,\text{in}\,\,
R^{n+1}\pi'_*(p_1^*\Omega_{\mathcal{C}/U}\otimes
p_2^*\Omega_{\mathcal{Y}_U}^{n})
\end{eqnarray}
where $g(\gamma)$ is the image of $\gamma\in
R^n\pi'_*\Omega_{\mathcal{C}\times_U\mathcal{Y}_U}^n$ in
$R^n\pi'_*(p_1^*\Omega_{\mathcal{C}/U}\otimes
p_2^*\Omega_{\mathcal{Y}_U}^{n-1})$ and
$\widetilde{\delta[\mathcal{D}]}$ is any lift of $[\mathcal{D}]$ in
$H^1(\rho^*\Omega_U)=\mathcal{H}^{0,1}_\mathcal{C}\otimes \Omega_U$.

We now observe that for $g(\gamma)\in
R^n\pi'_*(p_1^*\Omega_{\mathcal{C}/U}\otimes
p_2^*L^{s-1}\Omega_{\mathcal{Y}_U}^{n-1}$, and
$\widetilde{\delta[\mathcal{D}]}\in H^1(\rho^*\Omega_U)$ we have
\begin{eqnarray}
\label{eq8avrilfin} I(\widetilde{\delta[\mathcal{D}]}\cdot
g(\gamma))=w(<\delta[\mathcal{D}],\gamma_1>)\,\,\text{in}\,\,
R^n\pi_*\Omega_{\mathcal{Y}_U}^n/L^{s+1}R^n\pi_*\Omega_{\mathcal{Y}_U}^n.
\end{eqnarray}
In order to prove  (\ref{eq8avrilfin}), we choose  a class $\mu\in
H^1(\Omega_{\mathcal{C}/U})$ which does not vanish in
$H^0(U,R^1\rho_*\Omega_{\mathcal{C}/U})$ (for example the first
Chern class of a relatively ample line bundle on $\mathcal{C}$). The
class $g(\gamma)\in R^n\pi'_*(p_1^*\Omega_{\mathcal{C}/U}\otimes
p_2^*L^{s-1}\Omega_{\mathcal{Y}_U}^{n-1})$ can be written as
$$g(\gamma)=p_1^*\mu\cdot p_2^*\eta+p_2^*\gamma_1,
$$
with $\eta\in R^{n-1}\pi_*( L^{s-1}\Omega_{\mathcal{Y}_U}^{n-1})$ and $\gamma_1\in \mathcal{H}^{1,0}_\mathcal{C}\otimes
R^n\pi_*( L^{s-1}\Omega_{\mathcal{Y}_U}^{n-1})$.
We then find
that
$$I(\widetilde{\delta[\mathcal{D}]}\cdot p_1^*\mu\cdot p_2^*\eta)=\rho_*(\widetilde{\delta[\mathcal{D}]}\cdot\mu)\cdot \eta
$$
which is $0$ since we assumed $U$ affine so $\widetilde{\delta[\mathcal{D}]}\cdot\mu=0$
in $H^2(\mathcal{C},\Omega_{\mathcal{C}/U})$. On the other hand, it is clear
that for
$$\widetilde{\delta[\mathcal{D}]}\in \Omega_U\otimes H^1(\mathcal{O}_\mathcal{\mathcal{C}})=\Omega_U\otimes\mathcal{H}^{0,1}_\mathcal{C}\,
$$
we have
$$ I(\widetilde{\delta[\mathcal{D}]}\cdot p_2^*\gamma_1)=w(<\widetilde{\delta[\mathcal{D}]},\gamma_1>)$$
because for any
$u\in R^1\rho_*\mathcal{O}_\mathcal{C}$ and $\omega\in \mathcal{H}^{1,0}_\mathcal{C}$, with cup-product
$u\cdot\omega\in H^1(\Omega_{\mathcal{C}/U})$,
we have $I(u\cdot\omega)=<u,\omega>$.

Combining (\ref{eq8avrilprod}) and (\ref{eq8avrilfin}) clearly proves the desired formula
(\ref{formdesiree8avril}).

  \cqfd

We now turn to a general family $\mathcal{Y}\rightarrow B$. The proof works
in fact exactly in the same way, except that we have to replace the
natural splitting of the Leray filtration on  $H^n(Y_t,{\Omega_{\mathcal{Y}}^n}_{\mid Y_t})$
 given by the character decomposition by another one, which is not canonical, and depends only on the
 choice of a polarization on the family $\mathcal{Y}\rightarrow B$, that is a line bundle
 $\mathcal{L}$ on $\mathcal{Y}$ which is ample on the fibers $\mathcal{Y}_t$.
 \begin{prop} \label{decompleray} Assume the base $B$ is smooth quasi-projective. Given a relative polarization
 $\mathcal{L}$ on $\mathcal{Y}\stackrel{\pi}{\rightarrow}B$, there is a
 canonically induced splitting of the filtration on
 $H^q(\mathcal{Y}_t, {\Omega_\mathcal{Y}^p}_{\mid \mathcal{Y}_t})$ or of the sheaves $R^q\pi_*{\Omega_\mathcal{Y}^p}$
  induced by the filtration $L$ on
 ${\Omega_\mathcal{Y}^p}_{\mid Y_t}$. This splitting
 is functorial with respect to pull-back maps.
 \end{prop}

The proof of proposition \ref{decompleray} is given below.

 Choosing the pulled-back line bundle $p_2^*\mathcal{L}$ on the pulled-back family
 $\mathcal{C}\times_U\mathcal{Y}$, we get  a similar decomposition
 of $H^n(\mathcal{C}_t\times \mathcal{Y}_t,{\Omega_{\mathcal{C}\times_U\mathcal{Y}}^n}_{\mid \mathcal{C}_t\times \mathcal{Y}_t})$
which splits the  filtration $L$ on this space, relative to the projection
$\mathcal{C}\times_U\mathcal{Y}\rightarrow \mathcal{C}$.

 The decompositions are
compatible with pull-back and push-forward maps, and the rest of the
argument works the same way as before. This concludes the proof of
Proposition \ref{theoinfi}.

\cqfd
{\bf Proof of proposition \ref{decompleray}}. We could do it
at hand, but it is quickier to see  this as a consequence of  Deligne's canonical splitting result  in
\cite{delseattle}.
Deligne proves that the choice of a relative  polarization on
$\pi:\mathcal{Y}\rightarrow B$ induces a canonical quasi-isomorphism
$$ R\pi_*\mathbb{Q}\cong \oplus_i R^i\pi_*\mathbb{Q}.$$
It follows that there are canonical cohomology classes
$\tilde{\pi}_k\in H^{2d}(\mathcal{Y}\times_B
\mathcal{Y},\mathbb{Q}),\,d=\text{dim}\,\mathcal{Y}_t$ whose images
in $H^0(B,R^{2d}(\pi,\pi)_*\mathbb{Q})$ are the K\"unneth projectors
$\pi_k$:

$$R^{2d}(\pi,\pi)_*\mathbb{Q}=\oplus _{i+j=2d}R^i\pi_*\mathbb{Q}\otimes R^j\pi_*\mathbb{Q}=$$
$$
\oplus _{j}Hom\,(R^j\pi_*\mathbb{Q},R^j\pi_*\mathbb{Q})\ni
\pi_j:=Id_{R^j\pi_*\mathbb{Q}}.$$ From the construction of
\cite{delseattle}, it is clear that these classes $\tilde{\pi}_k\in
H^{2d}(\mathcal{Y}\times_B \mathcal{Y},\mathbb{Q})$ are Hodge
classes, which means that they extend to any smooth compactification
of $\mathcal{Y}\times_B \mathcal{Y}$ and are Hodge classes there. It
then follows that these classes have Dolbeault counterparts:
$$\tilde{\pi}_k^{d,d}\in H^{d}(\mathcal{Y}\times_B \mathcal{Y},\Omega_{\mathcal{Y}\times_B \mathcal{Y}}^d).
$$
We can now use the $\tilde{\pi}_k^{d,d}$ as acting on the relative
Dolbeault cohomology of $\mathcal{Y}$. Namely, denoting by
$pr_1,\,pr_2$ the projections from $\mathcal{Y}\times_B \mathcal{Y}$
to $B$, we have maps:
$$ pr_{2*}\circ (\tilde{\pi}_k)\cup \circ  pr_1^*: R^q\pi_*\Omega_\mathcal{Y}^p\rightarrow R^q\pi_*\Omega_\mathcal{Y}^p.$$
They similarly act on the fibers at $t$:
$$ pr_{2*}\circ (\tilde{\pi}_k)\cup \circ  pr_1^*: H^q(\mathcal{Y}_t,{\Omega_\mathcal{Y}^p}_{\mid \mathcal{Y}_t})\rightarrow
 H^q(\mathcal{Y}_t,{\Omega_\mathcal{Y}^p}_{\mid \mathcal{Y}_t}).$$
These maps are compatible with the $L$-filtration and induce a
morphism of spectral sequences associated to the filtration $L$. On
the other hand, their action on $E_1^{s,q-s}=\Omega_{B,t}^s\otimes
H^{q-s}(\mathcal{Y}_t,{\Omega_\mathcal{Y}^p}_{\mid \mathcal{Y}_t})$
is induced by $\pi_k$, hence it is equal to $0$ for $p+q-s\not=k$,
and the identity if $p+q-s=k$, and similarly for
$E_2^{s,q-s}=E_\infty^{s,q-s}$. It follows that we have induced maps
$$pr_{2*}\circ (\tilde{\pi}_{p+q-s})\cup \circ  pr_1^*: L^sH^q(\mathcal{Y}_t,{\Omega_\mathcal{Y}^p}_{\mid \mathcal{Y}_t})\rightarrow
L^sH^q(\mathcal{Y}_t,{\Omega_\mathcal{Y}^p}_{\mid \mathcal{Y}_t})$$
which have the property that they induce after passing to the graded
$Gr_L$ the zero map on $E_2^{s',q-s'}$ for $s'>s$, and the identity
map on $E_2^{s,q-s}=E_\infty^{s,q-s}$. But then a sufficiently high
power of these maps vanishes on
$L^{s+1}H^q(\mathcal{Y}_t,{\Omega_\mathcal{Y}^p}_{\mid
\mathcal{Y}_t})$ and thus factors through $E_\infty^{s,q-s}$ to give
a retraction $$E_\infty^{s,q-s}\rightarrow
L^sH^q(\mathcal{Y}_t,{\Omega_\mathcal{Y}^p}_{\mid \mathcal{Y}_t}),$$
hence a canonical splitting
$$L^sH^q(\mathcal{Y}_t,{\Omega_\mathcal{Y}^p}_{\mid \mathcal{Y}_t})=L^{s+1}H^q(\mathcal{Y}_t,{\Omega_\mathcal{Y}^p}_{\mid \mathcal{Y}_t})
\oplus E_\infty^{s,q-s},\,E_\infty^{s,q-s}=E_2^{s,q-s}.$$

\cqfd
\section{ $1$-cycles modulo algebraic equivalence in Jacobians}
\subsection{Constructing test elements in the dual space}
We consider now a family $\mathcal{C}\rightarrow B$ of curves of
genus $g$, and the associated Jacobian fibration
$\pi:\mathcal{J}\rightarrow B$. We choose an embedding
$\mathcal{C}\subset \mathcal{J}$, which provides a codimension $g-1$
cycle $\mathcal{Z} \in\text{CH}^{g-1}(\mathcal{J})/alg$ and we want
to study the cycles $\mathcal{Z}_s\in
\text{CH}^{g-1}(\mathcal{J})_s/alg$. Our goal is to exhibit such
families with a non trivial infinitesimal invariant
$$\delta[\mathcal{Z}_s]_{alg}\in I_{s,g-1-s,g-1}:=
\frac{\text{Ker}\,
(\overline{\nabla}_{s,g-1-s,g-1})}{\Omega_B\wedge\text{Ker}\,(\overline{\nabla}_{s-1,g-s-1,g-1})}$$
for $s\geq 2$. Here $\overline{\nabla}$ is the
$\overline{\nabla}$-map for the Hodge bundles of the family
$\mathcal{J}\rightarrow B$. The infinitesimal variation of Hodge
 structure on the cohomology
 of  the Jacobians  $\mathcal{J}_t,\,t\in B$, is deduced from the one on the cohomology  of the curves $\mathcal{C}_t,\,t\in B$,
in the obvious way, namely it coincides with the latter in degree $1$,
and with its $k$-th exterior power
 for higher $k$.

The map
$$\overline{\nabla}_{s,p,q}:\Omega_B^s\otimes \mathcal{H}^{p,q}\rightarrow
\Omega_B^{s+1}\otimes \mathcal{H}^{p-1,q+1}$$ has for transposed map
$$^t\overline{\nabla}_{s+1,g-p+1,g-q-1}:\bigwedge^{s+1}T_B\otimes \mathcal{H}^{g-p+1,g-q-1}
\rightarrow \bigwedge^{s}T_B\otimes \mathcal{H}^{g-p,g-q}.$$ This
map is given by the formula
$$^t\overline{\nabla}_{s+1,g-p+1,g-q-1}
(u_1\wedge\ldots\wedge u_{s+1}\otimes \omega_1\wedge\ldots\wedge\omega_{g-p+1}\otimes \eta_1\wedge\ldots\wedge\eta_{g-q-1}$$
$$=\sum_{1\leq i\leq s+1,\,1\leq j\leq g-p+1}(-1)^{i+j}u_1\wedge\ldots\wedge \hat{u_i}\wedge\ldots\wedge u_{s+1}\otimes \omega_1
\wedge\ldots\wedge \hat{\omega_j}\wedge\ldots\wedge \omega_{g-p+1}$$
$$\otimes ^t\overline{\nabla}_{1,1,0}(u_i\otimes \omega_j)\wedge
\eta_1\wedge\ldots\wedge\eta_{g-q-1}.$$

Fix a point $b\in B$.
The dual of the space of infinitesimal invariants
$$\frac{\text{Ker}\,
(\overline{\nabla}_{s,g-1-s,g-1,b}:\Omega_{B,b}^s\otimes
\mathcal{H}^{g-1-s,g-1}_b\rightarrow
 \Omega_{B,b}^{s+1}\otimes \mathcal{H}^{g-2-s,g}_b )}{\Omega_{B,b}\wedge\text{Ker}\,
 (\overline{\nabla}_{s-1,g-s-1,g-1,b}):\Omega_{B,b}^{s-1}
 \otimes \mathcal{H}^{g-1-s,g-1}_b\rightarrow \Omega_{B,b}^{s}\otimes \mathcal{H}^{g-2-s,g}_b}$$
 is the quotient of
the space
\begin{eqnarray}
\label{eqkernel8avril}
I^*_{s,g-1-s,g-1}:=\text{Ker}\,(\bigwedge^sT_{B,b}\otimes
\mathcal{H}^{s+1,1}_b\rightarrow T_{B,b}\otimes
\frac{\bigwedge^{s-1}T_{B,b}\otimes
\mathcal{H}^{s+1,1}_b}{\text{Im}\,(^t\overline{\nabla}_{s,s+2,0,b})})
\end{eqnarray}
by the subspace
$\text{Im}\,(^t\overline{\nabla}_{s+1,s+2,0,b}:\bigwedge^{s+1}T_{B,b}\otimes
\mathcal{H}^{s+2,0}_b\rightarrow\bigwedge^sT_{B,b}\otimes
\mathcal{H}^{s+1,1}_b)$,
 where $ \mathcal{H}^{p,q}_b
=\bigwedge^pH^0(C,K_C)\otimes \bigwedge^q H^1(C,\mathcal{O}_C)$. Let
us describe the simplest possible elements in the space $
I^*_{s,g-1-s,g-1}$:
\begin{lemm}\label{ledata}
Assume we have
\begin{enumerate}
\item a $s$-dimensional subspace
$W=<u_1,\ldots,u_s>\subset T_{B,b}$,
\item a $s+1$-dimensional subspace $K=<\omega_1,\ldots,\omega_{s+1}>\subset \mathcal{H}^{1,0}_b$ and a
$s$-dimensional subspace $K_1=<\mu_1,\ldots,\mu_s>\subset \mathcal{H}^{1,0}_b$ ,
\item an element $\eta\in \mathcal{H}^{0,1}_b$,
\end{enumerate}
such that
\begin{enumerate}
\item For any $i\in \{1,\ldots,s\}$ and $j\in \{1,\ldots,s+1\}$, $^t\overline{\nabla}(u_i\otimes \omega_j)=0$ in
$\mathcal{H}^{0,1}_b$.
\item For any $i\in \{1,\ldots,s\}$ and $j\in \{1,\ldots,s\}$, $^t\overline{\nabla}(u_i\otimes \mu_j)=\lambda_{ij}\eta$ in
$\mathcal{H}^{0,1}_b$, where the $(s,s)$-matrix $(\lambda_{ij})$ is
invertible.

\end{enumerate}
Then $w:=u_1\wedge\ldots\wedge u_s\otimes
\omega_1\wedge\ldots\wedge\omega_{s+1}\otimes\eta$ belongs to the
subspace $I^*_{s,g-1-s,g-1}$ of (\ref{eqkernel8avril}).
\end{lemm}
{\bf Proof.} Indeed, we may  assume that $\lambda_{ij}=\delta_{ij}$,
changing the basis  $u_i$ if necessary. Then the image of
$u_1\wedge\ldots\wedge u_s\otimes
\omega_1\wedge\ldots\wedge\omega_{s+1}\otimes\eta$ in
$T_{B,b}\otimes \bigwedge^{s-1}T_{B,b}\otimes \mathcal{H}^{s+1,1}$
is equal to
$\sum_i(-1)^iu_1\wedge\ldots\wedge\hat{u_i}\wedge\ldots\wedge
u_s\otimes \omega_1\wedge\ldots\wedge\omega_{s+1}\otimes\eta$. On
the other hand, we have for any $i$:
$$u_1\wedge\ldots\wedge\hat{u_i}\wedge\ldots\wedge u_s\otimes \omega_1\wedge\ldots\wedge\omega_{s+1}\otimes\eta
$$ $$=(-1)^i\,{^t\overline{\nabla}_{s,s+2,0,b}} (u_1\wedge\ldots\wedge
u_{s}\otimes \mu_i\wedge \omega_1\wedge\ldots \wedge\omega_{s+1}).$$

\cqfd
\subsection{The Ikeda family}
The Ikeda  family \cite{ikeda} is simply the family of plane curves
which are cyclic covers of $\mathbb{P}^1$ of degree $d$ ramified
along a degree $d$ divisor. In other words, the general equation
takes the form (in homogeneous coordinates $Y,X_0,X_1$ on
$\mathbb{P}^2$):
$$Y^d=f(X_0,X_1)$$
for some homogeneous polynomial $f$ of degree $d$.
This family is parameterized by the quotient $U/PGl(2)$ where
$U\subset \mathbb{P}(H^0(\mathbb{P}^1,\mathcal{O}_{\mathbb{P}^1}(d)))$ is the open set parameterizing reduced
divisors, hence smooth plane curves.

The tangent space to this quotient at a general point $f$ is the quotient
$$R^d_f=S^d_X/J_f^d,$$
where $S^k_X:=H^0(\mathbb{P}^1,\mathcal{O}_{\mathbb{P}^1}(k))$ and
$J_f^k\subset S^k$ is the degree $k$ part of the Jacobian ideal of
$f$, generated by the partial derivatives $\frac{\partial
f}{\partial X_i}$.

We consider the universal family $\mathcal{C}\rightarrow B$, where
$B\subset U$ is any slice for the action of $PGl(2)$ and the
corresponding infinitesimal variation of Hodge structure at a point
$b\in B$ corresponding to an equation $f$. Let us exhibit for any
such $f$ data as in Lemma \ref{ledata} with $s=d-3$.

We recall from \cite{cagri} (see also \cite[6.2.1]{voisinbook}) that
the infinitesimal variation of Hodge structure of an hypersurface
defined by an equation $F$ is governed by the product in the
Jacobian ring $R_F$ of $F$, quotient of the ring
$\mathbb{C}[Y,X_0,X_1]$ by the partial derivatives $\frac{\partial
F}{\partial X_i}$ and $\frac{\partial f}{\partial Y}$. In the case
of a curve $C$ of degree $d$ in $\mathbb{P}^2$ defined by an
equation $F$, we get
$$H^{1,0}(C)\cong R_F^{d-3},\,\,H^{0,1}(C)\cong R_F^{2d-3}$$
and the infinitesimal variation of Hodge structure restricted to the family of plane deformations of $C$
identifies
$$^t\overline{\nabla}_{1,1,0}:T_{B,b}\otimes H^{1,0}(C)\rightarrow H^{0,1}(C)$$
to the multiplication map
$$R^d_F\otimes R^{d-3}_F\rightarrow R^{2d-3}_F,$$
where $B$ is a slice for the $PGl(3)$-action and $b\in B$ is the parameter for $C$.

Let now $F=Y^d-f(X_0,X_1)$.  Note that $\frac{\partial F}{\partial
X_i}=-\frac{\partial f}{\partial X_i}$, so that $R_F$ is a module
over $R_f$.

Consider the following subspaces: \begin{eqnarray}
\label{eqK11avril} W:=R^d_{f}\subset R^d_F,\,K:=R_f^{d-3}=S_X^{d-3}
\end{eqnarray}
and $K_1:=YR_f^{d-4}=YS^{d-4}_X$. Note that $\text{dim}\,W=d-3$,
$\text{dim}\,K=d-2$, $\text{dim}\,K_1=d-3$. We have
\begin{lemm} \label{leikinf} (i) For any
$u\in W$, $\omega\in K$, we have $^t\overline{\nabla}(u\otimes
\omega)=0$.

(ii) There is an element $\eta\in H^{0,1}(C)$, such that the map
$W\otimes K_1\rightarrow H^{0,1}(C)$, $u\otimes \omega \mapsto
^t\overline{\nabla}(u\otimes \omega)$ takes the form
$$^t\overline{\nabla}(u\otimes \omega)=\lambda(u,\omega)\eta,$$
where the bilinear form $\lambda$ gives a perfect pairing between
$W$ and $K_1$.
\end{lemm}
{\bf Proof.} (i) This follows from the fact that $R^{2d-3}_f=0$.

 (ii) We
choose  for $\eta$ the generator of the $1$-dimensional vector space
$YR^{2d-4}_f$, which means geometrically (cf. \cite{cagri}) that
\begin{eqnarray} \label{eqeta11avril}
 \eta=\text{Res}_C\,\frac{PY\Omega}{f^2}\,\,\text{in}\,\,H^{0,1}(C),
\end{eqnarray}
 with $\Omega=X_0 dX_1\wedge dY-X_1dX_0\wedge dY+YdX_0\wedge dX_1\in
 H^0(\mathbb{P}^2,K_{\mathbb{P}^2}(3))$.

Then (ii) follows from Macaulay's theorem (cf.
\cite[6.2.2]{voisinbook}), which says that $R^{2d-4}_f=\mathbb{C}$
and the multiplication map $R^d_f\otimes R^{d-4}_f\rightarrow
R^{2d-4}_f$ is a perfect pairing. \cqfd We get now the following
corollary. Here we consider the Ikeda family $\mathcal{C}\rightarrow
B$ (defined over a slice of the quotient map
$\mathbb{P}(S^d_X)\rightarrow \mathbb{P}(S^d_X)/PGl(2)$), the
associated Jacobian fibration $\pi:\mathcal{J}\rightarrow B$, an
embedding $\mathcal{C}\rightarrow \mathcal{J}$ giving rise to a
codimension $g-1$ cycle $\mathcal{Z}$ of $\mathcal{J}$ and its
(relative) Beauville component $\mathcal{Z}_{d-3}$.
\begin{coro}\label{coroikinf}  Choose a basis $u_1,\ldots,u_{d-3}$ of $W$ and a basis
$\omega_1,\ldots,\omega_{d-2}$ of $K$. Let $\eta$ be a generator of
$YR^{2d-4}_f=\mathbb{C}\subset H^1(C,\mathcal{O}_C)$, where the
curve $C$ is defined by $Y^d=f(X_0,X_1)$. Then
\begin{eqnarray}\label{w11avril}
w:=u_1\wedge\ldots\wedge u_{d-3}\otimes\omega_1\wedge\ldots\wedge
\omega_{d-2}\otimes\eta \end{eqnarray} belongs to  the  space
$I_{d-3,g-d+2,g-1}^*$, $g=\frac{d(d-3)}{2}+1$.

\end{coro}
{\bf Proof.} We only have to observe that the space $W\subset
H^1(C,T_C)$ is exactly the tangent space to the Ikeda family, so
that $w$ belongs to $\bigwedge^{d-3}T_{B,b}\otimes H^{d-2,1}(JC)$.
 Applying Lemmas \ref{ledata} and Lemma \ref{leikinf}, we also get that
 the image of
 $w$ in $T_{B,b}\otimes \bigwedge^{d-4}T_{B,b}\otimes
 H^{d-2,1}(JC)$ vanishes in the quotient
 $$T_{B,b}\otimes \frac{\bigwedge^{d-4}T_{B,b}\otimes
 H^{d-2,1}(JC)}{\text{Im}\,^t\overline{\nabla}_{d-3,d-1,0}},$$
 which proves the result.
 \cqfd
\subsection{End of the proof of Theorem \ref{main}}
This section is devoted to the proof of the following result:
\begin{prop}\label{proinfnonzero} Let $f\in S^d_X$ define a smooth divisor in $\mathbb{P}^1$.
Let $C$ be the curve defined by $Y^d=f(X_0,X_1)$ and let $w\in
I_{s,g-1-s,g-1}^*$ be the element constructed in Corollary
\ref{coroikinf}. Then,  the pairing
$$<\delta[\mathcal{Z}_{d-3}]_{b,alg},w>$$
is nonzero, where $b\in B$ is the point of the Ikeda family
$\mathcal{C}\rightarrow B$ parameterizing $C$.

\end{prop}
Note that Proposition \ref{proinfnonzero} implies Theorem
\ref{main}. Indeed, it implies that for the Ikeda  family
$\mathcal{C}\rightarrow B$, the infinitesimal invariant
$\delta[\mathcal{Z}_{d-3}]_{b,alg}$ of the corresponding codimension
$g-1$ cycle $\mathcal{Z}$ of $\mathcal{J}\rightarrow B$ is nonzero
at the general point of $B$, which implies by Theorem
\ref{theoinfiintro} that for the very general curve $C$
parameterized by $B$ (hence for the very general plane curve), the
Beauville component  $Z_{d-3}$  is non zero  in
$\text{CH}_1(JC)_{d-3}/alg$.

The proof of Proposition \ref{proinfnonzero} will use the following
description of the pairing $<\delta[\mathcal{Z}_{d-3}]_{b,alg},w>$.
Let $\mathcal{C}\subset \mathcal{J}\rightarrow B$ be a family of
curves of genus $g$ and let $w\in I_{s,g-1-s,g-1}^*$ be an  element
constructed as in Lemma \ref{ledata}. Assume for simplicity that
$\text{dim}\,B=s$ so that $W=T_{B,b}$, where the point $b$
parameterizes the curve $C:=\mathcal{C}_b$. Then the pairing
$$<\delta[\mathcal{Z}_s]_{b,alg},w>$$
is computed as follows : Since all the form $\omega_j\in H^0(C,K_C)$
satisfy $^t\overline{\nabla}(u_i\otimes \omega_j)=0$ in
$H^1(C,\mathcal{O}_C)$, they lift to sections of
$$H^0(C,\Omega_{\mathcal{C}\mid C})=H^0(C,\Omega_{\mathcal{J}\mid
C})=H^0(J,\Omega_{\mathcal{J}\mid J}),$$ where $J=JC=\mathcal{J}_b$.
For any $k$, the multiplication $\mu_k$ induces an endomorphism
$\mu_k^*$ of $H^0(J,\Omega_{\mathcal{J}\mid J})$ which is compatible
with the restriction map
$$H^0(J,\Omega_{\mathcal{J}\mid J})\rightarrow H^0(J,\Omega_{ J}),$$
the action of $\mu_k^*$ on the right hand term being the homothety
of factor $k$.
 It follows that there is an
unique lift $\tilde{\omega}_j$ of $\omega_j$ in
$H^0(J,\Omega_{\mathcal{J}\mid J})$ satisfying the property
$$\mu_k^*(\tilde{\omega}_j)=k\tilde{\omega}_j\,\,\text{for any}\,\,k.$$
The $d-2$-form
$$\tilde{\omega}_1\wedge\ldots\wedge\tilde{\omega}_{s+1}\in H^0(C,\Omega_{\mathcal{J}\mid
C}^{s+1})$$ restricts to a form $\mu\in
H^0(C,\Omega_{\mathcal{C}\mid C}^{s+1})$. There is an
isomorphism
$$\Omega_{\mathcal{C}\mid C}^{s+1}\cong K_C$$
which becomes canonical once $\bigwedge^sW$ is trivialized and we
trivialize $\bigwedge^sW$ using the multivector
$u_1\wedge\ldots\wedge u_s$. We thus we get from $w$ an element
$\mu\cdot\eta\in H^1(C,K_C)$.

\begin{lemm}\label{ledescription}
We have
$$<\delta[\mathcal{Z}_s]_{b,alg},w>=\int_C\mu\cdot\eta.$$
\end{lemm}
{\bf Proof.} In a general situation $\mathcal{C}\hookrightarrow
\mathcal{J}\rightarrow B$, the class $[\mathcal{Z}]_b^{g-1,g-1}\in
H^{g-1}(JC,\Omega^{g-1}_{\mathcal{J}\mid JC})$ is the Dolbeault
cohomology class of $\mathcal{C}\subset \mathcal{J}$ and the class
$[\mathcal{Z}]_{s,alg,b}$ is its component lying in the space
$$H^{g-1}(JC,\Omega^{g-1}_{\mathcal{J}\mid JC})_{s}:=\{\alpha\in H^{g-1}(JC,\Omega^{g-1}_{\mathcal{J}\mid JC}),\,\mu_k^*\alpha=k^{2g-2-s},\,\forall k\}.$$
For any class $\beta$ in the dual space
$$H^1(JC,\Omega^{N+1}_{\mathcal{J}\mid JC}),\,N:=\text{dim}\,B$$
we have
\begin{eqnarray}
 \label{eqforint12}
<[\mathcal{C}]_b^{g-1,g-1},\beta>=\int_C\beta_{\mid \mathcal{C}},
\end{eqnarray}
where $\beta_{\mid \mathcal{C}}$ is the image of $\beta$ in
$H^1(C,\Omega^{N+1}_{\mathcal{C}\mid C})$ which is identified to
$H^1(C,K_C)$ via a trivialization of $\bigwedge^N\Omega_{B,b}$.
Coming back to our situation above, so $N=s$, the statement thus
follows from the following facts:

1) The choice of lift $\tilde{\omega}_j\in
H^0(\Omega_{\mathcal{J}\mid JC})_1$ gives the unique lift
$$\tilde{w}:=\tilde{\omega}_1\wedge\ldots\wedge\tilde{\omega}_{s+1}\otimes\eta\in
H^1(JC,\Omega_{\mathcal{J}\mid JC})_{s+2}$$
 of $w\in Gr_L^{s+2}H^1(JC,\Omega_{\mathcal{J}\mid JC}^{s+1})$.

 2) Writing
 $[\mathcal{Z}]^{g-1,g-1}=\sum_i[\mathcal{Z}_i]^{g-1,g-1}$, then
 we have
 \begin{eqnarray}
 \label{eqcoupling12}
 <[\mathcal{Z}]_b^{g-1,g-1},\tilde{w}>=<[\mathcal{Z}_s]_b^{g-1,g-1},\tilde{w}>,
 \end{eqnarray}
 since the coupling is multiplied by $k^{2g}$ under $\mu_k^*$, so
 that
 all the other terms
 $<[\mathcal{Z}_{s'}]^{g-1,g-1},\tilde{w}>$ vanish.

The term on the right in (\ref{eqcoupling12}) gives a coupling
$<\delta[\mathcal{Z}_s]_b^{g-1,g-1},{w}>$ which, by the fact that
$w$ is in the kernel  (\ref{eqkernel8avril}), is equal to the
desired coupling $<\delta[\mathcal{Z}_{s}]_{b,alg}^{g-1,g-1},{w}>$.
On the other hand,  the term on the left is computed by formula
(\ref{eqforint12}). The equality (\ref{eqcoupling12}) thus proves
Lemma \ref{ledescription}.
 \cqfd

 {\bf Proof of Proposition \ref{proinfnonzero}.} We follow the description given in Lemma \ref{ledescription}.
 As
$^t\overline{\nabla}(u_i\otimes\omega_j)=0$ for all $i,\,j$, each
$\omega_j$ has a canonical lift $\tilde{\omega}_j\in
H^0(J,\Omega_{\mathcal{J}\mid J})$. This lift is determined by the
condition that $\mu_n^*\tilde{\omega}_j=k\tilde{\omega}_j$ for all
$k$.

We observe that the cyclic group $G:=\mathbb{Z}/d\mathbb{Z}$ with
generator $g$ identified to a $d$-th  primitive root of unity
$\zeta$ acts over $B$  on the families $\mathcal{C}$ and
$\mathcal{J}$, in a compatible way if we embed $\mathcal{C}$ in
$\mathcal{J}$ using one of the fixed points of the action of $G$ on
the curves. We claim that the forms $\tilde{\omega}_j$ satisfy
\begin{eqnarray}
\label{eq15avril}g^*\tilde{\omega}_j=\zeta \tilde{\omega}_j.
\end{eqnarray}
To see this, we write the exact sequence $$0\rightarrow
\Omega_{B,b}\otimes \mathcal{O}_C\rightarrow \Omega_{\mathcal{C}\mid
C}\rightarrow \Omega_C\rightarrow 0,$$ which induces the exact
sequence
$$0\rightarrow \Omega_{B,f}\rightarrow H^0(C,\Omega_{\mathcal{J}\mid
C})\rightarrow K\rightarrow 0$$ with $K\subset H^{1,0}(C)$ defined
as in (\ref{eqK11avril}). Then the claim follows from the fact that
the $\omega_j$, which are the images of $\tilde{\omega}_j$ in
$H^0(J,\Omega_{\mathcal{J}/B\mid J})=H^0(J,\Omega_{ J})=H^{1,0}(C)$
satisfy the property
$$g^*{\omega}_j= \zeta{\omega}_j,$$
and from the fact that the two actions $\mu_n^*$ and $g^*$ on the
space
$$H^0(J,\Omega_{\mathcal{J}\mid J})=H^0(C,\Omega_{\mathcal{J}\mid C})=H^0(C,\Omega_{\mathcal{C}\mid C}),$$
commute.

As an immediate consequence of (\ref{eq15avril}), we find that
$$g^*(\tilde{\omega}_1\wedge\ldots\wedge\tilde{\omega}_{d-2}) =\zeta^{d-2}
\tilde{\omega}_1\wedge\ldots\wedge\tilde{\omega}_{d-2}\,\,\text{in}\,\,H^0(\bigwedge^{d-2}\Omega_{\mathcal{J}\mid
C}),$$ hence a fortiori the restricted form
$$\mu:=\tilde{\omega}_1\wedge\ldots\wedge\tilde{\omega}_{d-2}\in
H^0(\bigwedge^{d-2}\Omega_{\mathcal{C}\mid C})=H^{1,0}(C)$$
satisfies
$$g^*\mu=\zeta^{d-2}\mu.$$

Next we observe that the form $\eta\in H^{0,1}(C)$ defined in
(\ref{eqeta11avril}) satisfies
$$g^*\eta=\zeta^2\eta$$
as follows immediately from the formula
$$\eta=\text{Res}_C\,\frac{PY\Omega}{f^2},$$
with $g^*Y=\zeta Y,\,g^*\Omega=\zeta\Omega$, $g^*f=f,\,g^*P=P$.

The two spaces $H^{1,0}(C)^{\zeta^{d-2}}$ and $H^{0,1}(C)^{\zeta^2}$
where $g^*$ acts respectively by multiplication by $\zeta^{d-2}$ and
$\zeta^2$ are dual and $1$-dimensional and the second space is
generated by $\eta$. It follows that the pairing
$<\delta[\mathcal{Z}_{d-3}]_{b,alg},w>$, which by Lemma
\ref{ledescription} is computed as the pairing
 $<\mu,\eta>$ between $\mu\in H^{1,0}(C)$ and $\eta\in H^{0,1}(C)$,
 vanishes if and only if
 $\mu$ vanishes identically.

 It thus suffices to prove that $\mu$ is nonzero.
This is done as follows: Recall that the bundle
$\Omega_{\mathcal{C}\mid C}$ has rank $d-2$ and possesses the $d-2$
sections $\tilde{\omega}_j$. Assume that $\mu$ vanishes in
$H^0(C,\bigwedge^{d-2}\Omega_{\mathcal{C}\mid C})=H^0(C,K_C)$. Then
these sections generate a subbundle $N\subset
\bigwedge^{d-2}\Omega_{\mathcal{C}\mid C}$ of rank $\leq d-3$, and
it follows that if $x\in C$ is a general point, there is a section
$\tilde{\omega}_x\in<\tilde{\omega}_1,\ldots,\tilde{\omega}_{d-2}>$
vanishing at $x$. As all the sections $\tilde{\omega}_j$ satisfy
$g^*\tilde{\omega}_j=\zeta\tilde{\omega}_j$, it follows that
$\tilde{\omega}_x$ also vanishes at $gx,g^2x,\ldots , g^{d-1}x$, so
that $\tilde{\omega}_x$ is in fact a section of
$\Omega_{\mathcal{C}\mid C}(-1)$. It is easy to check however using
the multiplication in the Jacobian ring $R_f$ that $
H^0(C,\Omega_{\mathcal{C}\mid C}(-1))=0$. Indeed, we have the exact
sequence $$0\rightarrow\Omega_{B,b}\otimes
\mathcal{O}_C(-1)\rightarrow \Omega_{\mathcal{C}\mid
C}(-1)\rightarrow K_C(-1)\rightarrow 0$$ and the induced map $$
H^0(C,K_C(-1))\rightarrow \Omega_{B,b}\otimes
H^1(C,\mathcal{O}_C(-1))=\text{Hom}\,(R^d_f,R^{2d-4}_F)$$ is given
by multiplication in the Jacobian ring of $F$.
 Such a section would
provide a nonzero element of
$$R^{d-4}_F
 =\oplus_{i=0}^{i=d-1}Y^iR^{d-4-i}_f$$
 which is  annihilated by multiplication
by all elements of $W=R^d_f$, and this does not exist by Macaulay's
Theorem.

 \cqfd
\section{Comparison with Ikeda's invariants}

\label{secdiscussion}

As already mentioned in the introduction, Ikeda \cite{ikeda} proves
the nonvanishing of the cycle $Z_{d-3}$ for the very general curve
in the Ikeda family of plane curves of degree $d$, modulo a certain
subgroup of $\text{CH}^{g-1}(JC)$. Here the cycle $Z_{d-3}$ is
obtained by  embedding of the curve $C$ in $JC$ using a general
point $p\in JC$.

The infinitesimal invariant used by Ikeda is the projection of
$$\delta[\mathcal{Z}_{d-3,b}]\in \frac{ \text{Ker}\overline{\nabla}_{d-3,g-d+2,g-1,b}\subset \Omega_{B,b}^{d-3}\otimes
H^{g-d+2,g-1}(JC)}{\text{Im}\,\overline{\nabla}_{d-4,g-d+3,g-2}}
$$

modulo the subspace
\begin{eqnarray}\label{eqikedasubspace}\Omega_{B,b}^{d-3}\otimes
\text{Ker}\,\overline{\nabla}_{0,g-d+2,g-1,b}\subset
\text{Ker}\overline{\nabla}_{d-3,g-d+2,g-1,b}.
\end{eqnarray}

In  contrast, our infinitesimal invariant
$\delta[\mathcal{Z}_{d-3,alg,b}]$ is the projection of
$$\delta[\mathcal{Z}_{d-3,b}]\in \frac{ \text{Ker}\overline{\nabla}_{d-3,g-d+2,g-1,b}\subset \Omega_{B,b}^{d-3}\otimes
H^{g-d+2,g-1}(JC)}{\text{Im}\,\overline{\nabla}_{d-4,g-d+3,g-2}}
$$

modulo the subspace \begin{eqnarray} \label{eqmyspace}
\Omega_{B,b}\wedge
\text{Ker}\,\overline{\nabla}_{d-4,g-d+2,g-1,b}\subset
\text{Ker}\overline{\nabla}_{d-3,g-d+2,g-1,b} \end{eqnarray}
 which
contains (\ref{eqikedasubspace}) by formula (\ref{eqnew5avril}).

If we consider (instead of the family of plane curves) the universal
family of curves $\mathcal{C}\rightarrow B$ of genus $g$, and
replace $d-3$ by any integer $s$ such that $g\geq 2s+1$, (this is
the range we are interested in by Conjecture \ref{conj}) then at the
general point $b\in B$, the subspace (\ref{eqikedasubspace}) can be
in fact explicitly computed, namely:
\begin{prop} \label{prokerikeda} Let $C=\mathcal{C}_b$ be a general curve of genus $g\geq 2s+1$.
Then the kernel
\begin{eqnarray}
\label{eqnabla12avril}\text{Ker}\,(\overline{\nabla}_{0,g-s-1,g-1,b}:H^{g-s-1,g-1}(JC)\rightarrow
\Omega_{B,b}\otimes H^{g-s-2,g})
\end{eqnarray}
 is equal to the image
$$\text{Im}\,([C]*:H^{g-s,g}(JC)\rightarrow H^{g-s-1,g-1}(JC)) ,$$
where $[C]*$ is  the Pontryagin product with the class of $C$ in
$H^{g-1,g-1}(JC)$.

\end{prop}

Notice that $\text{Im}\,([C]*:H^{g-s,g}(JC)\rightarrow
H^{g-s-1,g-1}(JC))$ is also equal to the nonprimitive part of the
cohomology
$$\text{Im}\,(\theta^{g-s-1}:H^{0,s}(JC)\rightarrow
H^{g-s-1,g-1}(JC)).$$

 {\bf Proof of Proposition \ref{prokerikeda}.} The proposition above is in  fact a
geometric translation of our result in \cite{voisinsyz}, namely the
vanishing of the Koszul cohomology group $ K_{k, 1}(C,K_C)$ for a
general curve of genus $g$ in the range $g \leq 2k+1$. This
vanishing says that in this range,
 the kernel
of the Koszul differential \begin{eqnarray} \label{eqkoszul}
\text{Ker}\, (\bigwedge^kH^0(C,K_C)\otimes H^0(C,K_C)\rightarrow
\bigwedge^{k-1}H^0(C,K_C)\otimes H^0(C,2K_C)) \end{eqnarray}
 is equal
to the image of the natural map (also a Koszul differential)
$$\bigwedge^{k+1}H^0(C,K_C)\rightarrow\bigwedge^kH^0(C,K_C)\otimes
H^0(C,K_C).$$ If we now identify $H^0(C,K_C)$ with
$\bigwedge^{g-1}H^{0,1}(C)= \bigwedge^{g-1}H^{0,1}(C)$, the Koszul
differential (\ref{eqkoszul}) identifies to the map
$\overline{\nabla}_{0,k,g-1},\,k=g-s-1$ of (\ref{eqnabla12avril})
because, as is well known, the multiplication map $$
H^0(C,K_C)\otimes H^0(C,K_C)\rightarrow H^0(C,2K_C)$$ seen in the
form
$$ H^0(C,K_C)\rightarrow \text{Hom}(
H^0(C,K_C),H^0(C,2K_C))=H^{0,1}(C)\otimes H^0(C,2K_C)$$ identifies
to the map $\overline{\nabla}_{0,1,0}$.

For $g\geq 2s+1$, we have $g\leq 2(g-s-1)+1$.
 It thus follows from the above that the
kernel of $\overline{\nabla}_{0,g-s-1,g-1}$ is equal to the image of
$\bigwedge^{g-s}H^0(C,K_C)$. We finally observe that the natural map
$$\bigwedge^{g-s}H^0(C,K_C)\rightarrow
\bigwedge^{g-s-1}H^0(C,K_C)\otimes H^0(C,K_C)$$ identifies to the
 map $[C]*:H^{g-s,g}(JC)\rightarrow H^{g-s-1,g-1}(JC)$
 given by Pontryagin product with the class of $C$.

 We refer to \cite[8.1.4]{apnagel} for a more detailed version of
this Hodge theoretic interpretation of the vanishing of canonical
syzygies.
 \cqfd
Proposition \ref{prokerikeda} shows that the Ikeda subspace
(\ref{eqikedasubspace}) is  too small to make the Ikeda
infinitesimal invariant of the  cycles $\mathcal{Z}_s\in
\text{CH}^{g-1}(\mathcal{J})$ for $g\geq 2s+1$ invariant under
relative translation of the cycle. Indeed, choose a degree $1$
codimension $g$ cycle $\Gamma$ in $\mathcal{J}$, and consider the
cycle
$$\mathcal{Z}_{\Gamma}:=\mathcal{Z}*\Gamma,$$
where $*$ is the relative Pontryagin product  over $B$ for cycles in
$\mathcal{J}$. Then
$$\mathcal{Z}_{\Gamma}-\mathcal{Z}=\mathcal{Z}*(\Gamma-0_B),$$
where $0_B$ is the zero section of $\mathcal{J}$, and since the
$0$-cycles $\Gamma_b-0_b\in \text{CH}_0(\mathcal{J}_b)$ are of
degree $0$, they are algebraically equivalent to $0$, and so are the
cycles $\mathcal{Z}_{\Gamma,b}-\mathcal{Z}_b$.

If we now introduce the relative Beauville decomposition
$$\Gamma=\sum_{i}\Gamma_i, \,\Gamma_0=0_B$$ of $\Gamma$, we find that
$$\mathcal{Z}_{\Gamma,s}=\sum_{i\leq s}\mathcal{Z}_{s-i}*\Gamma_i.$$
In this sum, the first term is $\mathcal{Z}_s$ and the last term is
$\mathcal{Z}_{0}*\Gamma_s$, and its restriction to $\mathcal{J}_b$
over the point $b\in B$ belongs by definition to the subgroup
$Z_0F^s\text{CH}_{g-1}(\mathcal{J}_b)$. So its infinitesimal
invariant belongs by Ikeda's results to the Ikeda subspace
(\ref{eqikedasubspace}).
 The other terms however have
infinitesimal invariants which can be shown, using Proposition
\ref{prokerikeda} and Ceresa's Theorem, not to belong to the Ikeda
subspace (\ref{eqikedasubspace}) for $g\geq 2s+1$, (while of course
they belong to the subspace (\ref{eqmyspace}) by Proposition
\ref{theoinfi}, since $\mathcal{Z}_{s-i}*\Gamma_i$ is algebraically
equivalent to $0$ for $i>0$).

\end{document}